\documentclass[11pt,a4paper]{article}

\usepackage{amsmath,amssymb, bbm}
\usepackage{color}
\newcommand{\R}{\mathbb{R}}
\newcommand{\N}{\mathbb{N}}
\newcommand{\Z}{\mathbb{Z}}
\newcommand{\A}{\mathbb{A}}
\newcommand{\B}{\mathbb{B}}
\def\cA{{\mathcal A}}
\def\cC{{\mathcal C}}
\def\cD{{\mathcal D}}
\def\cS{{\mathcal S}}
\def\cF{{\mathcal F}}
\newcommand{\ee}{\varepsilon}
\renewcommand{\aa}{\alpha}
\renewcommand{\div}{{\rm div}\,}

\newcommand{\Sum}{\displaystyle \sum}

\def\d{\partial}
\def\ddl{\dot \Delta_l}
\def\ddj{\dot \Delta_j}
\def\ddq{\dot \Delta_q}

\def\tilde{\widetilde}
\def\hat{\widehat}
\newcommand{\D}{\Delta}
\newcommand{\DF}{\Delta_F}
\newcommand{\La}{\Lambda}
\newcommand{\n}{\nabla}
\newcommand{\G}{\Gamma}
\newcommand{\Ge}{G^e}
\newcommand{\Fe}{F^e}
\newcommand{\Om}{\Omega}

\newcommand{\GL}{\Gamma_L}

\newcommand{\fr}{\frac{1}{r}}
\newcommand{\frb}{\frac{2}{r}}

\newcommand{\uj}{u_j}
\newcommand{\ujt}{\tilde{u_j}}

\newcommand{\Fj}{F_j^e}
\newcommand{\Gj}{G_j^e}
\newcommand{\Rj}{R_j}

\newcommand{\Sjt}{\tilde{S_j}}
\newcommand{\Rjt}{\tilde{R_j}}
\newcommand{\SjtL}{\tilde{S_j^L}}
\newcommand{\SjtNL}{\tilde{S_j^{NL}}}
\newcommand{\Fjt}{\tilde{F_j^e}}
\newcommand{\Gjt}{\tilde{G_j^e}}

\newcommand{\psj}{\psi_{j,t}}
\newcommand{\psjm}{\psi_{j,t}^{-1}}
\newcommand{\psjto}{\psi_{j,\tau}}
\newcommand{\no}{\nu_0}

\newcommand{\Phie}{\Phi_\ee}
\newcommand{\ve}{v_\ee}
\newcommand{\Ue}{U_\ee}

\newcommand{\Thee}{\theta_\ee}

\newtheorem{thm}{Theorem}
\newtheorem{lem}{Lemma}

\newtheorem{prop}{Proposition}
\newtheorem{defi}{Definition}

\newtheorem{rem}{Remark}

\title{A priori estimates for the 3D quasi-geostrophic system}

\author{Fr\'ed\'eric Charve\footnote{Universit\'e Paris-Est Cr\'eteil, Laboratoire d'Analyse et de Math\'ematiques Appliqu\'ees (UMR 8050), 61 Avenue du G\'en\'eral de Gaulle, 94 010 Cr\'eteil Cedex (France). E-mail: frederic.charve@u-pec.fr}}

\date{}
\begin{document}

\maketitle

\begin{abstract} The present article is devoted to the 3D dissipative quasi-geostrophic system ($QG$). This system can be obtained as limit model of the Primitive Equations in the asymptotics of strong rotation and stratification, and involves a non-radial, non-local, homogeneous pseudo-differential operator of order 2 denoted by $\G$ (and whose semigroup kernel reaches negative values). After a refined study of the non-local part of $\G$, we prove apriori estimates (in the general $L^p$ setting) for the 3D $QG$-model. The main difficulty of this article is to study the commutator of $\G$ with a Lagrangian change of variable. An important application of these a priori estimates, providing bound from below to the lifespan of the solutions of the Primitive Equations for ill-prepared blowing-up initial data, can be found in a companion paper.
\end{abstract}

\section{Introduction}

\subsection{Presentation of the model}

A very rich litterature is devoted to the $2D$ quasi-geostrophic system:
\begin{equation}
\begin{cases}
\d_t \theta +v\cdot \n \theta +|D|^\aa \theta=0,\\
\theta_{|t=0}=\theta^0,
\end{cases}
\label{QG2}
\tag{$2DQG$}
\end{equation}
where $\theta$ represents the scalar potential temperature, $\aa\in[0,2]$, and the two-dimensional velocity $v$ is determined from $\theta$ by:
$$
v=(-\d_2|D|^{-1}, \d_2|D|^{-1}) \theta= (-R_2, R_1) \theta.
$$
The operators $R_i$ ($i=1,2$) are Riesz transforms. Due to the fractional diffusion term $|D|^\aa \theta$, three cases have to be distinguished. First, the sub-critical case $(\aa>1)$, which is now well known: global existence and uniqueness for arbitrary initial data are stated in various spaces. In the critical case $(\aa=1)$, results are proved for regular initial data with smallness assumptions. The super-critical case $(\aa<1)$ seems to be harder and more recently studied. Non-exhaustively, we refer to \cite{Chae, CLM, Const, Cordoba, TH2, TH3, Ju, LMM, KNV, LMM, M}.
\\

Pretty much less is done for the $3D$ quasi-geostrophic system on which we will focus in this article. As for the $2D$ version it is related to geophysical fluids. The system of Primitive Equations, that describes a geophysical fluid which is located at the surface of the Earth, writes as follows:
\begin{equation}
\begin{cases}
\d_t \Ue +\ve\cdot \n \Ue -L \Ue +\frac{1}{\ee} \cA \Ue=\frac{1}{\ee} (-\n \Phie, 0),\\
\div \ve=0,\\
{\Ue}_{|t=0}=U_{0,\ee}.
\end{cases}
\label{PE}
\tag{$PE_\ee$}
\end{equation}
The unknows are $\Ue =(\ve, \Thee)=(\ve^1, \ve^2, \ve^3, \Thee)$ ($\ve$ is the velocity of the fluid and $\Thee$ denotes the scalar potential temperature, which is related to the density fluctuation) and $\Phie$ is called the geopotential and includes the centrifugal force and the pressure. The diffusion operator $L$ is given by
$$
L\Ue \overset{\mbox{def}}{=} (\nu \D \ve, \nu' \D \Thee),
$$
where $\nu, \nu'>0$ are the kinematic viscosity and the thermal diffusivity. The antisymmetric matrix $\cA$ is defined by
$$
\cA \overset{\mbox{def}}{=}\left(
\begin{array}{llll}
0 & -1 & 0 & 0\\
1 & 0 & 0 & 0\\
0 & 0 & 0 & F^{-1}\\
0 & 0 & -F^{-1} & 0
\end{array}
\right).
$$
We also assume that the initial data converges towards $U_0=(v_0, \theta_0)= \underset{\ee \rightarrow 0}{\mbox{lim }} U_{0,\ee}$.
\\

The small parameter $\ee$ measures the importance of the Coriolis force (modelized by the first $3\times3$ diagonal submatrix of $\cA$, which corresponds to the vector product $\ee^{-1} v_\ee \times e_3$ of the velocity with the third unit vector $e_3$) and of the vertical stratification of the density induced by the gravity (corresponding to the other terms in the matrix and the parameter $F\in ]0,1]$). We chosed here the small parameter so that these two concurrent phenomena are of the same importance.

We will not give more details about this system and adress the reader to \cite{Sadourny, Pedlosky, BMN5, FC5, FCpochesLp} for a more precise presentation of the physical models. And we also refer to \cite{CDGG, CDGG2, CDGG3, Dutrifoy2, IG1} concerning the rotating fluids system and to \cite{TH4, TH5} for similar methods for the Boussinesq system.

The limit system when the small parameter $\ee$ goes to zero, is called the 3D quasi-geostrophic system and consists in a transport-diffusion equation coupled with a Biot-Savart-type law. This system writes as follows (we refer to \cite{Chemin2, FC2, FC, FC3, FC4, FC5, FCpochesLp, Dutrifoy1} for studies of the asymptotics as $\ee$ goes to zero):
\begin{equation}
\begin{cases}
\d_t \Om +v.\n \Om -\G \Om =0\\
U=(v, \theta)=(-\partial_2,\partial_1,0,-F\partial_3)\DF^{-1} \Om,\\
\Om_{|t=0}=\Om_0
\end{cases}
\tag{QG}
\label{QG1}
\end{equation}
where the operator $\G$ is defined by:
$$
\G \overset{def}{=} \D \DF^{-1} (\nu \d_1^2 +\nu \d_2^2+ \nu' F^2 \d_3^2),
$$
with $\DF=\d_1^2 +\d_2^2 +F^2 \d_3^2$, and $\Om=\d_1 U^2 -\d_2 U^1 -F \d_3 U^4 =\d_1 v^2 -\d_2 v^1 -F \d_3 \theta$ and $\Om_0=\d_1 v_0^2 -\d_2 v_0^1 -F \d_3 \theta_0$ where $U_0=(v_0, \theta_0)$ is the limit as $\ee$ goes to zero of the initial data $U_{0,\ee}$.

\begin{rem}
\sl{Except in the cases $\nu=\nu'$ (where $\G=\nu \D$, see \cite{FC3} and \cite{Dutrifoy2} for example) or $F=1$ (where $\G= \nu \d_1^2 +\nu \d_2^2+ \nu' \d_3^2$, we refer to \cite{Chemin2}) the operator $\G$ is a non-local diffusion operator of order 2.
}
\end{rem}

\begin{rem}
\sl{$(QG)$ can be rewritten into a system close to the classical incompressible Navier-Stokes system, satisfied by the velocity $U=(-\partial_2,\partial_1,0,-F\partial_3)\DF^{-1} \Om$ (see \cite{FC2, FCpochesLp}).}
\end{rem}
We refer to \cite{FC} where we obtain existence of global Leray solutions if the initial velocity $U_0 =(-\partial_2,\partial_1,0,-F\partial_3)\DF^{-1} \Om_0 \in L^2$. Similarly the Fujita-Kato theorem is easily adapted and we have local existence of a strong solution if $U_{0, QG} \in \dot{H}^\frac{1}{2}$, and global existence for small data.
\begin{rem}
\sl{We refer to \cite{FC2} where we proved that in fact, using the quasigeostrophic structure, if $U_{0,QG} \in H^1$ then we have a global strong solution without assuming any smallness condition for the initial data. We emphasize that it is a very remarkable fact, which relates the $3D$-QG system to the $2D$ incompressible Navier-Stokes system more than to the $3D$ version. Another connection between these systems is that as for the vorticity in the $2D$-case, the $3D$-QG system has no stretching term $\Om\cdot \n v$.
}
\end{rem}
As in \cite{FC2, FC3} the long-term aim is to obtain qualitative results for the solutions of the Primitive Equations. Due to the remarkable properties of the limit system \eqref{QG1}, we are able to obtain similar properties for the solutions of the Primitive Equations provided that $\ee>0$ is small enough. For example in \cite{FC2}, we used the fact that if the quasi-geostrophic part of the initial data $U_0$ (independant of $\ee$ in the cited paper) is in $H^1$, the limit system \eqref{QG1} has a unique global solution and the same is true for the Primitive Equations if $\ee>0$ is small enough (without any smallness conditions on the initial data). In other words, the fast rotation and strong stratification help the Primitive Equations, which are very close to the $3D$-Navier-Stokes system, to have global strong solutions. Put it differently, adding a penalized skew-symmetric term to the $3D$ incompressible Navier-Stokes system ($3DNS$) helps filtering the fast oscillations. This stabilizes the system and allows it to have global solutions without smallness conditions (which is not possible with these assumptions for $(3DNS)$ where we are not able to prove global existence of strong solutions if the initial data is $H^1$ and large).
\\

As explained in the present article we will obtain new a priori estimates for System \eqref{QG1} in the general case (no relations between the kinematic viscosity $\nu$ and the thermal diffusivity $\nu'$). These estimates will be the key to obtain a bound from below for the lifespan $T_\ee^*$ of the solutions of the Primitive Equations for ill-prepared initial data, with large and possibly blowing-up (in $\ee$) norms in $\dot{H}^\frac{1}{2}$ (that is assumptions where the Fujita-Kato theorem cannot provide us informations on the Lifespan). More precisely, we will show that $T_\ee^* \geq \gamma \ln(\ln|\ln \ee|)$. This work is done in our companion paper \cite{FCpochesLp} and then generalizes the result from \cite{FC3} where we made the simplifying assumption $\nu=\nu'$ (turning the operator $\G$ into $\nu \D$). The Prandtl number, which can be defined as $Pr=\nu/\nu'$, can take values far from $1$ so that it is important to overcome the restriction from \cite{FC3} and treat the case where there is no relation between $\nu$ and $\nu'$. 

We will also use these a priori estimates in forth-coming works on the $3DQG$-system.

\subsection{Statement of the main results}

We consider the following transport-diffusion system:

\begin{equation}
\begin{cases}
\d_t u +v.\n u -\G u =\Fe,\\
u_{|t=0}=u_0
\end{cases}
\label{TDQG}
\end{equation}
Let us introduce $M_{visc}=\frac{\max(\nu, \nu')}{\min(\nu, \nu')}$.

\begin{prop}\sl{($L^p$-estimates) Assume that $u$ solves \eqref{TDQG} on $[0,T]$ with $u_0 \in L^p$ and that $\|v\|_{L_T^\infty L^6} \leq C'$ (for some constant $C'$) with $\div v=0$. Then there exists a constant $D$ (depending on $F$, $M_{visc}$ and $C'$) such that for all $t\in [0,T]$,
\begin{equation}
\|u\|_{L_t^\infty L^p} \leq D^t (\|u_0\|_{L^p}+\int_0^t \|\Fe(\tau)\|_{L^p} d\tau).
\end{equation}  
}
\label{estimLp}
\end{prop}

\begin{thm}\sl{(Smoothing effect) Assume that $u$ solves \eqref{TDQG} on $[T_1, T_2]$ with $v$ satisfying $\div v=0$ and $\|v\|_{L^\infty({[T_1,T_2]}, L^6)} \leq C'$, $u(T_1)\in L^p$, $\Fe \in L_{loc}^1 L^p$ (for $p\in[1,\infty]$). There exist two constants $C$ and $C_F$ such that if $T_2-T_1>0$ is so small that:
\begin{enumerate}
\item $2CC'(T_2-T_1)^{\frac{1}{4}} \leq \no^3$,
\item $e^{C\int_{T_1}^{T_2} \|\nabla v (\tau)\|_{L^\infty}}-1 \leq \frac{1}{C_F M_{visc}}$.
\end{enumerate}
Then, for all $r\in[1,\infty]$, there exists a constant $C_{r,F}>0$ such that for all $t\in [T_1, T_2]$,
\begin{equation}
(\no r)^{\fr} \|u\|_{\tilde{L}^r([T_1,t], B_{p, \infty}^\frb)} \leq C_{r, F} \left(\|u(T_1)\|_{L^p} +\int_{T_1}^t \|\Fe(\tau)\|_{L^p} d\tau \right).
\end{equation}
}
\label{thC2}
\end{thm}

\begin{rem}
\sl{In the particular case $r=1, p=\infty$ we obtain that:
\begin{equation}
\no \|u\|_{\tilde{L}^1([T_1,t], C_*^2)} \leq C_F \left(\|u(T_1)\|_{L^\infty} +\int_{T_1}^t \|\Fe(\tau)\|_{L^\infty} d\tau \right).
\end{equation}
}
\end{rem}

\begin{thm}\sl{(a priori estimates)
Let $s\in ]-1,1[$. Assume that $u$ solves \eqref{TDQG} on $[T_1, T_2]$ with $v$ satisfying $\div v=0$ and $\|v\|_{L^\infty({[T_1,T_2]}, L^6)} \leq C'$, $u(T_1)\in B_{p,\infty}^s$. Assume in addition that the external force term can be decomposed into $\Fe+\Ge$, with $\Fe \in \tilde{L}^1([T_1,T_2], B_{p,\infty}^s)$ (for $p\in[1,\infty]$) and $\Ge \in \tilde{L}^\infty([T_1,T_2], B_{p,\infty}^{s+\frac{2}{r}-2})$ for $r\in[1,\infty]$ with $s+\frac{2}{r}\in]-1,1[$. There exist two constants $C_s$ and $C_F$ such that if $T_2-T_1>0$ is so small that:
\begin{enumerate}
\item $2CC'(T_2-T_1)^{\frac{1}{4}} \leq \no^3$,
\item $e^{C\int_{T_1}^{T_2} \|\nabla v (\tau)\|_{L^\infty}}-1 \leq \frac{1}{C_F M_{visc}}$,
\item $T_2-T_1 +\int_{T_1}^{T_2} \|\n v\|_{L^\infty} d\tau \leq C_s$
\end{enumerate}
Then, there exists a constant $C_{\no,F}>0$ such that for all $t\in [T_1, T_2]$,
\begin{multline}
(\no r)^\frac{1}{r} \|u\|_{\tilde{L}^r([T_1,t], B_{p, \infty}^{s+\frac{2}{r}})}\\
\leq C_{\no, F} \left(\|u(T_1)\|_{B_{p,\infty}^s} +\|\Fe\|_{\tilde{L}^1([T_1,t],B_{p,\infty}^s)} +\frac{1}{\no}\|\Ge\|_{\tilde{L}^\infty ([T_1,t],B_{p,\infty}^{s+\frac{2}{r}-2})} \right).
\end{multline}
}
\label{thCs}
\end{thm}

\begin{rem}
\sl{The time globalization can be done as T. Hmidi did in \cite{TH1} and only introduces a multiplicative factor $e^{C(t-T_1+\int_{T_1}^t \|\n v(\tau)\| d\tau)}$ in the results. We choose to write on a small time interval because this is the form we will use in \cite{FCpochesLp}.}
\end{rem}

\begin{rem}
\sl{In the energy case $p=2$, these estimates are easier to obtain.}
\end{rem}

\begin{rem}
\sl{These estimates have been obtained by T. Hmidi in \cite{TH1} for the Navier-Stokes system and adapted in \cite{FC3} for the Primitive Equations in the case $\nu=\nu'$ where $\G=\nu \D$. In the general case treated by the present article, all the difficulties result from the fact that $\G$ is a non local linear operator. The main features will be to obtain estimates for commutators involving $\G$ or to cope with terms of the form $\G(uv)$.}
\end{rem}

The article is structured as follows: the second section is devoted to the proof of Proposition \ref{estimLp} ($L^p$ estimates). In Section $3$ we rewrite the non-local operator $\G$ into a more practical form and obtain various properties and commutator estimates that we use in Sections $4$ and $5$ to prove the theorems. In the appendix we gathered an introduction to the Besov spaces, general properties for flows, and a study of particular diffeomorphisms introduced in Section $3$.

\begin{rem}
\sl{In this paper we will denote as $C$ a universal constant and as $C_{F,s}$ (for example) a constant which only depends on $F$ and $s$. Even if from line to line this constant may change (multilplication by another constant...) we will still denote it as $C_{F,s}$.}
\end{rem}

\section{$L^p$ estimates}

The aim of this section is to prove Proposition \ref{estimLp}.

\subsection{Basic semigroup estimates}

The object of this preliminary section is to prove the following results on the semigroup generated by the operator $\G$. Basically $\G$ is "close to" some $\aa \D$ and the associated semigroup is expected to be "close to" the classical heat semigroup $e^{\aa t\D}$ (see \ref{decpgamma} for more details about $\aa$). In fact even when $\nu \sim\nu'$, $e^{t\hat{\G}(\xi)}$ is close to the gaussian function $e^{-t\aa |\xi|^2}$ but its Fourier inverse reaches negative values and has a $L^1$-norm strictly larger than $1$. This is why we cannot expect a genuine maximum principle for $\G$ and get the exponential factor $D^t$ in Proposition \ref{estimLp}. For the same reason we were unable to use the Trotter formula or the arguments from \cite{Cordoba, TH5}.

\begin{prop}
\sl{There exists a constant $C=C_{F,visc}>0$ depending on $F$ and $M_{visc}=\frac{\max(\nu, \nu')}{\min(\nu, \nu')}$ such that for all $p\in[1,\infty]$ and $u\in L^p$, we have:
$$
\|e^{t\G} u\|_{L^p} \leq C \| u\|_{L^p}.
$$
}
\label{semig1}
\end{prop}

\begin{prop}
\sl{There exists a constant $C=C_{F,visc, c_0, C_0}>0$ depending on $F$, $M_{visc}$, $c_0$ and $C_0$ such that for all $p\in[1,\infty]$, $\lambda>0$ and $u\in \cS'$ with $\mbox{supp }\hat{u} \in \lambda \cC$, where $\cC$ is the annulus $\cC(0, c_0, C_0)$, we have:
$$
\|e^{t\G} u\|_{L^p} \leq C e^{-\frac{c_0^2}{8} \no t \lambda^2} \| u\|_{L^p}.
$$
}
\label{semig2}
\end{prop}

\textbf{Proof of Proposition \ref{semig1}:} we will follow the lines of the classical proofs for the heat equation case (we refer for example to \cite{Dbook}), except that we must be cautious when applying the integration by parts argument due to the presence of a rational function as we have:
$$
\hat{\G u}(\xi)= -\frac{|\xi|^2}{|\xi|_F^2} (\nu \xi_1^2+ \nu \xi_2^2 +\nu' F^2 \xi_3^2) \hat{u}(\xi)\overset{\mbox{def}}{=} -q(\xi)\hat{u}(\xi).
$$
This implies that $e^{t\G} u= K_t * u$, where we define the kernel $K_t(x)$ for all $t,x$ by:
$$
K_t(x)=\cF^{-1} (e^{tq(\xi)})(x) =\frac{1}{(2\pi)^3} \int_{\R^3} e^{ix\cdot \xi} e^{-tq(\xi)} d\xi.
$$
Let us introduce $M=\frac{\nu}{\no}$ and $M'=\frac{\nu'}{\no}$ with $\no=\min(\nu, \nu')$, then we can write that:
$$
\begin{cases}
\min(M, M')=1,\\
\max(M,M')=M_{visc}=\frac{\max(\nu, \nu')}{\min(\nu, \nu')},
\end{cases}
$$
and for all $t,x$,
$$
K_t(x)=\frac{1}{(2\pi)^3} \int_{\R^3} e^{ix\cdot \xi} e^{-q_0(\xi \sqrt{\no t})} d\xi, \quad \mbox{with} \quad q_0(\xi)= \frac{|\xi|^2}{|\xi|_F^2} (M \xi_1^2+ M \xi_2^2 +M' F^2 \xi_3^2).
$$
Performing the change of variable $\xi=\frac{\eta}{\sqrt{\no t}}$, we obtain that for all $t,x$,
$$
K_t(x)=\frac{1}{\sqrt{\no t}^3} K_1(\frac{x}{\sqrt{\no t}}), \quad \mbox{with} \quad K_1(x)=\frac{1}{(2\pi)^3} \int_{\R^3} e^{ix\cdot \xi} e^{-q_0(\xi)} d\xi.
$$
Then as $e^{t\G} u= K_t * u = \frac{1}{\sqrt{\no t}^3} K_1(\frac{.}{\sqrt{\no t}}) * u$, and thanks to convolution estimates,
$$
\|e^{t\G} u\|_{L^p} \leq \|K_1\|_{L^1} \|u\|_{L^p},
$$
so that Proposition \ref{semig1} will be proved as soon as we can bound $\|K_1\|_{L^1}$. To do this, as in \cite{Dbook} we write:
\begin{multline}
K_1(x)=C\frac{1}{(1+|x|^2)^2} \int_{\R^3} (I_d-\D_{\xi})^2 (e^{ix\cdot \xi}) e^{-q_0(\xi)} d\xi\\
=C\frac{1}{(1+|x|^2)^2} \int_{\R^3} e^{ix\cdot \xi} (I_d-\D_{\xi})^2 (e^{-q_0(\xi)}) d\xi.
\end{multline}
\begin{rem}
\sl{Note that in \cite{Dbook} the operator $I_d-\D_{\xi}$ is applied $d$ times (in the general $\R^d$ case). In fact, it is sufficient to apply it $[d/2]+1$ times. However, as $q_0$ is a rational fraction($q$ is a polynomial in the heat case), deriving once more would lead to non integrable terms (at zero)}
\end{rem}
Collecting the polynomial part of $q_0$ gives:
$$
q_0(\xi)=\bigg(M \xi_1^2+ M \xi_2^2 +\big((1-F^2)M+F^2 M'\big) \xi_3^2\bigg) -(M-M')F^2(1-F^2)\frac{\xi_3^4}{|\xi|_F^2},
$$
From this, we compute the derivatives of $q_0$ and we can write that for all $\xi\neq 0$:
$$
q_0(\xi)=Q_2(\xi), \quad \n q_0(\xi)=Q_1(\xi),\quad \mbox{and } \D q_0(\xi)=Q_0(\xi),
$$
where $Q_i$ denotes a rational homogeneous fraction in $\xi$ of degree $i$ (the denominator is a power of $|\xi|_F$) and whose coefficients depend on $F$ and $M_{visc}$. This allows us to obtain that (with the same notations):
$$
\begin{cases}
(I_d-\D)e^{-q_0}=(1+Q_0+Q_2)e^{-q_0},\\
(I_d-\D)^2 e^{-q_0}=(Q_{-2}+Q_0+Q_2+Q_4)e^{-q_0},
\end{cases}
$$
then, using the fact that $-q_0(\xi)\leq -\min(M,M') |\xi|^2=-|\xi|^2$ we can estimate $K_1(x)$ for all $x\in \R^3$ by:
$$
|K_1(x)|\leq \frac{C_{F,M_{visc}}}{(1+|x|^2)^2} \int_{\R^3} (|\xi|^{-2} +1 +|\xi|^2 +|\xi|^4) e^{-|\xi|^2} d\xi \leq \frac{C_{F,M_{visc}}}{(1+|x|^2)^2},
$$
which implies that $\|K_1\|_{L^1}\leq C_{F,M_{visc}}$. $\blacksquare$
\begin{rem}
\sl{We can similarly prove that for all $k\in \N$, there exists a constant $C_{F,M_{visc},k}$ such that for all $x\in \R^3$:
$$
|\n^k K_1(x)|\leq \frac{C_{F,M_{visc},k}}{(1+|x|^2)^2},
$$
so that $\n^k K_1\in L^p$ for all $p\in [1,\infty]$.
}
\label{rqnabla}
\end{rem}
\textbf{Proof of Proposition \ref{semig2}: } as above, we use the same arguments as in \cite{Dbook}. As $\mbox{supp }\hat{u} \subset \lambda \cC$, if $\phi$ is a smooth function compactly supported in a bigger annulus $\cC'=\cC(0, c_0/2,2 C_0)$ such that $\phi\equiv 1$ on $\cC$, then we have $e^{t\G} u= g_t * u$, where for all $t,x$:
$$
g_t(x)=\cF^{-1} (e^{t\no q_0(.)} \phi(\frac{.}{\lambda}))(x) =\frac{1}{(2\pi)^3} \int_{\R^3} e^{ix\cdot \xi} e^{-t\no q_0(\xi)} \phi(\frac{\xi}{\lambda}) d\xi =\lambda^3 h_{\no t \lambda^2} (\lambda x),
$$
with
$$
h_{\tau}(x)= \frac{1}{(2\pi)^3} \int_{\cC'} e^{i x\cdot \xi} \phi(\xi) e^{-\tau q_0(\xi)} d\xi.
$$
Using the same integration by part as before, we write that:
$$
h_{\tau}(x)= \frac{C}{(1+|x|^2)^2} \int_{\cC'} e^{i x\cdot \xi} (I_d-\D)^2 (\phi(\xi) e^{-\tau q_0(\xi)}) d\xi.
$$
The rest of the proof is classical, computing $(I_d-\D)^2 \big(\phi(\xi) e^{-\tau q_0(\xi)}\big)$, using the fact that $\xi \in \cC'$ and that $-q_0(\xi)\leq -|\xi|^2$, we obtain that for all $x\in \R^3$,
$$
|h_\tau (x)| \leq \frac{C_{F,M_{visc}, c_0, C_0}}{(1+|x|^2)^2} \int_{\cC} (1+\tau^4) e^{-\tau \frac{c_0^2}{4}} d\xi \leq \frac{C_{F,M_{visc}, c_0, C_0}}{(1+|x|^2)^2} e^{-\tau \frac{c_0^2}{8}},
$$
then we immediately get that:
$$
\|g_t\|_{L^1}= \|h_{\no t \lambda ^2}\|_{L^1} \leq C_{F,M_{visc}, c_0, C_0} e^{-\frac{c_0^2}{8} \no t  \lambda^2},
$$
which ends the proof. $\blacksquare$

\subsection{Proof of Proposition \ref{estimLp}}

As $u$ solves system \eqref{TDQG}, the Duhamel form gives that for all $t\in[0,T]$ (using that $v$ is divergence-free),
$$
u(t)= e^{t\G} u_0 +\int_0^t e^{(t-\tau)\G} \big(-\div(v\otimes u)(\tau)+\Fe(\tau)\big) d\tau.
$$
Thanks to Proposition \ref{semig1} we can estimate the $L^p$-norm, there exists a constant $C$ such that for all time:
\begin{multline}
\|u(t)\|_{L^p} \leq C\left( \|u_0\|_{L^p} +\int_0^t \|\Fe(\tau)\|_{L^p} d\tau\right)\\
+ \int_0^t \|\frac{1}{\sqrt{\no(t-\tau)}^3} K_1(\frac{.}{\sqrt{\no(t-\tau)}}) * \div(v\otimes u)\|_{L^p} d\tau.
\end{multline}
In the convolution term we have:
\begin{multline}
\|\frac{1}{\sqrt{\no(t-\tau)}^3} K_1(\frac{.}{\sqrt{\no(t-\tau)}}) * \div(v\otimes u)\|_{L^p}\\
\leq \|\frac{1}{\sqrt{\no(t-\tau)}^4} (\n K_1)(\frac{.}{\sqrt{\no(t-\tau)}}) * (v\otimes u)\|_{L^p}\\
\leq \frac{1}{\sqrt{\no(t-\tau)}^4} \|(\n K_1)(\frac{.}{\sqrt{\no(t-\tau)}})\|_{L^a}  \|v\otimes u\|_{L^b},
\end{multline}
thanks to the classical convolution estimates, if $a,b\in[1,\infty]$ satisfy $\frac{1}{a} +\frac{1}{b}=1+\frac{1}{p}$. And if $r$ satisfies: $\frac{1}{p} +\frac{1}{r}=\frac{1}{b}$, that is in fact we simply ask that $\frac{1}{a} +\frac{1}{r}=1$, then $\|v\otimes u\|_{L^b} \leq \|v\|_{L^r} \|u\|_{L^p}$. Taking $a=6/5$ and $r=6$ we get:
\begin{multline}
\|u(t)\|_{L^p} \leq C\left( \|u_0\|_{L^p} +\int_0^t \|\Fe(\tau)\|_{L^p} d\tau\right)\\
+\int_0^t \frac{1}{\sqrt{\no(t-\tau)}^4} \sqrt{\no(t-\tau)}^{3\cdot \frac{5}{6}} \|\n K_1\|_{L^\frac{6}{5}} \|v(\tau)\|_{L^6} \|u(\tau)\|_{L^p} d\tau.
\end{multline}
Thanks to Remark \ref{rqnabla}, the integral can be estimated by
$$
C(F, M_{visc}) \int_0^t \no^{-\frac{3}{4}}(t-\tau)^{-\frac{3}{4}} \|v\|_{L_t^\infty L^6} \|u\|_{L_t^\infty L^p} d\tau \leq C'' \no^{-\frac{3}{4}} t^{\frac{1}{4}} \|u\|_{L_t^\infty L^p},
$$
where the constant $C''$ depends on $C'$, $F$ and $M_{visc}$. Then
$$
\|u\|_{L_t^\infty L^p} \leq C\left( \|u_0\|_{L^p} +\int_0^t \|\Fe(\tau)\|_{L^p}d\tau \right) +C''\no^{-\frac{3}{4}} t^\frac{1}{4} \|u\|_{L_t^\infty L^p},
$$
and if $t$ is so small that $C''\no^{-\frac{3}{4}}t^\frac{1}{4} \leq \frac{1}{2}$, we obtain:
\begin{equation}
\|u\|_{L_t^\infty L^p} \leq C \left(\|u_0\|_{L^p}+\int_0^t \|\Fe(\tau)\|_{L^p} d\tau\right).
\label{Linftpetit}
\end{equation}
Finally when $t$ is large, we globalize this result thanks to a subdivision $0=T_0<T_1<...<T_N=t$ such that for all $i\in \{0,...,N-1\}$,
\begin{equation}
C''(T_{i+1}-T_i)^\frac{1}{4} \sim \frac{1}{2}\no^\frac{3}{4},
\label{subdiv1}
\end{equation}
then for all $i\in \{0,...,N-1\}$ and all $t'\in [T_i, T_{i+1}]$,
$$
\|u\|_{L_{[T_i,t']}^\infty L^p} \leq C \left(\|u(T_i)\|_{L^p}+\int_{T_i}^{t'} \|\Fe(\tau)\|_{L^p} d\tau\right),
$$
which classically implies that:
$$
\|u\|_{L_t^\infty L^p} \leq C^N \left(\|u_0\|_{L^p}+\int_0^t \|\Fe(\tau)\|_{L^p} d\tau\right).
$$
Thanks to \eqref{subdiv1}, $N\sim (2C'')^4\no^{-3} t$, we obtain the desired result with $D=C^\frac{(2C'')^4}{\no^3}$. $\blacksquare$

\section{Around the operator $\Gamma$}

The proof of Theorems \ref{thC2} and \ref{thCs} follows classical lines: as in the work of T. Hmidi for the Navier-Stokes system (see \cite{TH1}) or in \cite{CD, CHVP, Corder} for compressible models, we localise the equations thanks to the dyadic operators, then we perform a Lagrangian change of variable, localizing again we obtain the desired estimates. As in \cite{CHVP, Corder} all the difficulty lies in the study of the commutator of the nonlocal operator $\G$ with the Lagrangian change of variable. In what follows we will skip details on what is classical and focus on this point. 

\subsection{Rewriting of the operator: singular integrals}

In the previous studies of the Primitive Equations and the 3D-QG system, we dealt with the operator $\G$ in energy spaces (only using its symbol, see \cite{FC, FC2, FC4, FC5}) or in particular cases where $\G$ is reduced to the classical Laplacian ($\nu=\nu'$, see \cite{FC3}). In the present article we will need a more handy expression for this operator and the next section is devoted to rewrite $\G$ as a singular integral (in contrast with the cases of \cite{CHVP, Corder} where the kernels have nicer properties).

Let us begin by collecting the pure local and non-local parts out of $\G$, from its definition, we easily write that (recall that $\D_F =\d_1^2 +\d_2^2 +F^2 \d_3^2$):
\begin{multline}
\G= \D \DF^{-1} (\nu \d_1^2 +\nu \d_2^2+ \nu' F^2 \d_3^2)= \D \DF^{-1} (\nu \DF +(\nu'-\nu) F^2 \d_3^2)\\
=\D (\nu I_d +(\nu'-\nu) F^2 \d_3^2 \DF^{-1} ) =\nu \D +(\nu'-\nu) F^2 \d_3^2 \D\DF^{-1} .
\end{multline}
As $\D= \DF+(1-F^2) \d_3^2$, we obtain:
$$
\G= \nu \D +(\nu'-\nu) F^2 \d_3^2 (I_d+(1-F^2) \d_3^2\DF^{-1} )= \GL +(\nu-\nu') F^2(1-F^2) \La^2,
$$
where we define the following operators:
\begin{equation}
\begin{cases}
\GL= \nu \d_1^2 +\nu \d_2^2+ \left((1-F^2)\nu +F^2 \nu'\right) \d_3^2,\\
\La= \d_3^2 (-\DF)^{-\frac{1}{2}}.
\end{cases}
\label{decpgamma}
\end{equation}
\begin{rem}
\sl{Dealing with $\d_3^4 (-\DF)^{-1}$ would lead to an obstruction but (as in \cite{TH3}) simply studying the square root of this operator (which is the derivative of a Riesz operator) will provide the properties we need.
}
\end{rem}
Let us recall the following result:
\begin{prop} (1.29, see \cite{Dbook} p.23)
\sl{If $|.|$ denotes the canonic euclidean norm in $\R^d$, then for each $\sigma \in ]0,d[$, there exists a constant $C_{d,\sigma}>0$ such that:
$$
\cF (|.|^{-\sigma}) =C_{d,\sigma} |.|^{\sigma-d}.
$$
}
\end{prop}
We deduce from this that, with $d=3$ and $\sigma=2$:
$$
\cF^{-1}(\frac{1}{|\xi|_F})= \frac{C}{F} \frac{1}{x_1^2+ x_2^2 +\frac{1}{F^2}x_3^2} =  \frac{C}{F} \frac{1}{|x|_{\frac{1}{F}}^2},
$$
where $C>0$ denotes a universal constant and we introduce, for $x\in \R^3$ and $\aa\neq 0$,
\begin{equation}
|x|_\aa^2 =x_1^2 +x_2^2 + \aa^2 x_3^2.
\label{normodif}
\end{equation}
We can then define the distribution $T$ as follows: for all $\psi \in \cD$, 
$$
<T,\psi> =\int_{\R^3} \frac{\psi(x)}{|x|_{\frac{1}{F}}^2} dx =\int_{\R^3} \frac{\psi(x)}{x_1^2+ x_2^2 +\frac{1}{F^2}x_3^2} dx.
$$
As we want to compute $\La =\frac{C}{F} \d_3^2 T$, we begin with writing (thanks to the Lebesgue theorem): for all $\psi \in \cD$, 
$$
<T,\psi> =\lim_{\ee \rightarrow 0} \int_{|x|\geq \ee} \frac{\psi(x)}{x_1^2+ x_2^2 +\frac{1}{F^2}x_3^2} dx.
$$
Thanks to the Green formula and the fact that the surface integral goes to zero with $\ee$, we easily obtain that:
$$
<\d_3 T,\psi> =<PV\bigg( \d_3 (\frac{1}{|x|_{\frac{1}{F}}^2)})\bigg), \psi>=\lim_{\ee \rightarrow 0} \int_{|x|\geq \ee} -\frac{2}{F^2}\frac{x_3}{\big(x_1^2+ x_2^2 +\frac{1}{F^2}x_3^2\big)^2} \psi(x)dx.
$$
Then, deriving once more:
$$
<\d_3^2 T,\psi> =-<\d_3 T,\d_3\psi> =\lim_{\ee \rightarrow 0} \bigg(\int_{|x|\geq \ee}  K_0(x)\psi(x)dx + J_\ee\bigg),
$$
with the kernel (degree $-4$):
$$
K_0(x)=\d_3^2 (\frac{1}{|x|_{\frac{1}{F}}^2)}) =-\frac{2}{F^2}\frac{x_1^2+ x_2^2 -\frac{3}{F^2}x_3^2}{\big(x_1^2+ x_2^2 +\frac{1}{F^2}x_3^2\big)^3},
$$
and the surface integral:
\begin{multline}
J_\ee =\int_{|x|=\ee} -\frac{2}{F^2}\frac{x_3^2}{|x|\big(x_1^2+ x_2^2 +\frac{1}{F^2}x_3^2\big)^2} \psi(x) d\sigma_{\ee}(x)\\
=-\frac{2}{F^2} \int_0^{2\pi} \int_{-\frac{\pi}{2}}^\frac{\pi}{2} \frac{\sin^2 \varphi \cos \varphi}{\big(\cos^2 \varphi +\frac{1}{F^2}\sin^2 \varphi\big)^2} \frac{\psi(\ee \cos \varphi \cos \theta, \ee \cos \varphi \sin \theta, \ee \sin \varphi)}{\ee } d\varphi d\theta.
\end{multline}

\begin{lem}
\sl{The previous integral satisfies:
$$
I_\ee \overset{\mbox{def}}{=} J_\ee+\frac{4\pi}{F^2}\frac{1}{\ee} \int_{-\frac{\pi}{2}}^\frac{\pi}{2} \frac{\sin^2 \varphi \cos \varphi}{\big(\cos^2 \varphi +\frac{1}{F^2}\sin^2 \varphi\big)^2} \psi(0) d\varphi \underset{\ee \rightarrow 0}{\longrightarrow} 0.
$$}
\end{lem}
\textbf{Proof: } the quantity is equal to:
$$
I_\ee= -\frac{2}{F^2} \int_{-\frac{\pi}{2}}^\frac{\pi}{2} \frac{\sin^2 \varphi \cos \varphi}{\big(\cos^2 \varphi +\frac{1}{F^2}\sin^2 \varphi\big)^2} \frac{\psi(\ee \cos \varphi \cos \theta, \ee \cos \varphi \sin \theta, \ee \sin \varphi)-\psi(0)}{\ee } d\varphi,
$$
and thanks to the Lebesgue theorem,
\begin{multline}
\lim_{\ee \rightarrow 0} I_\ee =-\frac{2}{F^2} \n \psi(0) \cdot \int_0^{2\pi} \int_{-\frac{\pi}{2}}^\frac{\pi}{2} \frac{\sin^2 \varphi \cos \varphi}{\big(\cos^2 \varphi +\frac{1}{F^2}\sin^2 \varphi\big)^2}
\left(\begin{array}{l}
\cos \varphi \cos \theta\\
\cos \varphi \sin \theta\\
\sin \varphi
\end{array}
\right) d\varphi d\theta\\
=-\frac{2}{F^2} \n \psi(0) \cdot
\left(\begin{array}{l}
0\\
0\\
0
\end{array} \right)=0,
\end{multline}
as the three integrals are zero. $\blacksquare$

Moreover, thanks to the Green formula, for any fixed $\ee$, 
$$
-\frac{4\pi}{F^2}\frac{1}{\ee} \int_{-\frac{\pi}{2}}^\frac{\pi}{2} \frac{\sin^2 \varphi \cos \varphi}{\big(\cos^2 \varphi +\frac{1}{F^2}\sin^2 \varphi\big)^2} \psi(0) d\varphi =-\int_{|x|\geq \ee}  K_0(x) dx,
$$
so that we have, for all $\psi \in \cD$,
$$
<\d_3^2 T,\psi> =\lim_{\ee \rightarrow 0} \int_{|x|\geq \ee}  K_0(x)\big(\psi(x)-\psi(0)\big)dx,
$$
and we finally end up with the following expression of $\La$: for all function $f\in \cS$, and all $x\in \R^3$,
\begin{multline}
\La f (x)=\frac{C}{F} \big(\d_3^2 T *f \big)(x) =\frac{C}{F} <\d_3^2 T, f (x-.)>\\
=\lim_{\ee \rightarrow 0} \int_{|y|\geq \ee}  K(y)\big(f(x-y)-f(x)\big) dy,
\label{exprdiff1}
\end{multline}
with the kernel $K$ defined for all $y\in \R^3$ by:
\begin{equation}
K(y)= -\frac{2C}{F^3}\frac{y_1^2+ y_2^2 -\frac{3}{F^2}y_3^2}{\big(y_1^2+ y_2^2 +\frac{1}{F^2}y_3^2\big)^3}.
\label{Kexpr}
\end{equation}
Note that as the degree of $K$ is $-4$, the integral diverges at the origin, so that we still need to desingularize this expression. An easy way to do this is (as used in \cite{CHVP}) to perform the change of variable $y\mapsto -y$ and we obtain (as $K$ is even) that for all $x\in \R^3$,
\begin{multline}
\La f (x) =\lim_{\ee \rightarrow 0} \int_{|y|\geq \ee}  K(y)\big(f(x-y)-f(x)\big) dy =\lim_{\ee \rightarrow 0} \int_{|y|\geq \ee}  K(y)\big(f(x+y)-f(x)\big) dy\\
=\lim_{\ee \rightarrow 0} \int_{|y|\geq \ee}  \frac{K(y)}{2}\bigg(f(x-y)+f(x+y)-2f(x)\bigg) dy\\
=\frac{1}{2}\int_{\R^3}  K(y)\bigg(f(x-y)+f(x+y)-2f(x)\bigg) dy.
\label{exprdiff2}
\end{multline}
As in \cite{CHVP} we will use the following alternative characterizations of Besov norms (see \cite{Dbook}), if we denote $\tau_a:f\mapsto f(.-a)$:
\begin{thm}\sl{(\cite{Dbook}, $2.36$)
 Let $s \in ]0,1[$ and $p,r\in [1,\infty]$. There exists a constant $C$ such that for any $u\in \cC_h'$,
$$
C^{-1} \|u\|_{\dot{B}_{p,r}^s}\leq \|\frac{\|\tau_{-y}u -u\|_{L^p}}{|y|^s}\|_{L^r (\R^d; \frac{dy}{|y|^d})} \leq C \|u\|_{\dot{B}_{p,r}^s}.
$$}
\end{thm}
and when $s=1$, we have to use finite differences of order $2$ instead of order $1$:
\begin{thm}\sl{(\cite{Dbook}, $2.37$)
 Let $p,r\in [1,\infty]$. There exists a constant $C$ such that for any $u\in \cC_h'$,
$$
C^{-1} \|u\|_{\dot{B}_{p,r}^1}\leq \|\frac{\|\tau_{-y}u +\tau_y u -2u\|_{L^p}}{|y|}\|_{L^r (\R^d; \frac{dy}{|y|^d})} \leq C \|u\|_{\dot{B}_{p,r}^1}.
$$}
\label{thdiff2}
\end{thm}
Unsurprisingly, we use the latter to retrieve that for all $p\in [1,\infty]$, and all $f\in\dot{B}_{p,1}^1(\R^3)$, then $\La f \in L^p(\R^3)$ and
\begin{equation}
\|\La f\|_{L^p} \leq C_F \|f\|_{\dot{B}_{p,1}^1}.
\label{estimLaBesov}
\end{equation}
\begin{rem}
\sl{The main reason why we studied $\La$ instead of directly $\La^2$ is that, as the latter is a pseudo-differential operator of integer order 2, dealing with Besov spaces would require to obtain a similar expression but with finite differences of order 3, introducing heavy problems with singular integrals.
}
\end{rem}

\subsection{Localization}

The proofs of Theorems $1$ and $2$ share the same starting point. As announced, if $u$ solves System \eqref{TDQG}, we first apply the dyadic operator $\D_j$ (we refer to the appendix for the Littlewood-Paley theory), and we obtain, as in \cite{TH1, Dlagrangien, CD, CHVP, Corder} that $\uj \overset{\mbox{def}}{=} \D_j u$ satisfies:
\begin{equation}
\begin{cases}
\d_t \uj +S_{j-1}v.\n \uj -\G \uj =\Fj+\Gj+\Rj,\\
{\uj}_{|t=0}=\D_j u_0= u_{0,j},
\end{cases}
\end{equation}
where the well-known remainder term is defined as follows
\begin{equation}
\Rj=(S_{j-1}v -v).\n \uj +[v.\n, \D_j]u =S_{j-1}v.\n \uj -\D_j(v\cdot \n u).
\label{Rj}
\end{equation}
\begin{rem}
\sl{We emphasize that the terms $S_{j-1}v.\n \uj$ and $R_j$ have their frequencies located in a ring of size $2^j$ if $j\geq 1$, and as explained in the proofs,  we will use it only for $j\geq N_0$ with $N_0$ large enough.}
\end{rem}

The next step is as in the cited works to perform a Lagrangian change of variable. We first define $\psj$ the flow associated to the regularized advection velocity $S_{j-1}v$, that is the solution of the following system:
\begin{equation}
\begin{cases}
\d_t \psj (x)=S_{j-1}v(t,\psj (x)),\\
\psi_{j,0}(x)=x.
\end{cases}
\label{defpsi}
\end{equation}
Then, introducing $\ujt=\uj \circ \psj$, that is for all $t,x$, $\ujt(t,x)=\uj(t,\psj(x))$ (same notation for the other quantities), we get: 
\begin{equation}
\begin{cases}
\d_t \ujt -\G \ujt =\Fjt+\Gjt+\Rjt+\Sjt,\\
\ujt_{|t=0}=u_{0,j},
\end{cases}
\end{equation}
where we have defined $\Sjt$ as the following commutator:
\begin{equation}
\Sjt=(\G \uj)\circ \psj-\G(\uj \circ \psj).
\end{equation}
Decomposing $\G$ into its pure local and non-local parts (see \eqref{decpgamma}):
$$
\G=\GL +(\nu-\nu') F^2(1-F^2) \La^2,
$$
we decompose the commutator $\Sjt$ into $\Sjt=\SjtL+\SjtNL$, with:
\begin{equation}
\begin{cases}
\vspace{1mm}
\SjtL=(\GL \uj)\circ \psj-\GL(\uj \circ \psj),\\
\SjtNL=(\nu-\nu')F^2(1-F^2) \bigg((\La^2 \uj)\circ \psj-\La^2(\uj \circ \psj)\bigg).
\end{cases}
\end{equation}
Obviously, we got rid of the advection term, but ended with functions that are not frequency localized anymore, and the next step is to localize once again and taking advantage of the fact that (as $\div S_{j-1}v=0$, $\psj$ is volume-preserving):
\begin{multline}
\|u_j\|_{L^p} =\|\ujt \circ \psjm\|_{L^p} =\|\big(S_{j-N_0}\ujt\big) \circ \psjm +\big(\sum_{q\geq j-N_0} \D_q\ujt\big) \circ \psjm\|_{L^p}\\
\leq \|S_{j-N_0}\ujt\|_{L^p} +\sum_{q\geq j-N_0} \|\D_q\ujt\|_{L^p}.
\label{decpjl}
\end{multline}
We refer to the next section for the low frequencies, as they will be dealt using Lemma $2.6$ from \cite{Dbook} exactly as in \cite{TH1, CD, CHVP, Corder}, and we focus here on the high frequencies, for wich we apply the dyadic operator $\D_l$ to the previous system. We get that for all $l\geq j-N_0$:
\begin{equation}
\begin{cases}
\d_t \D_l\ujt -\G \D_l\ujt =\D_l \bigg(\Fjt+\Gjt+\Rjt+\SjtL+\SjtNL\bigg),\\
\D_l\ujt_{|t=0}=\D_l u_{0,j} =\D_l \D_j u_0,
\end{cases}
\end{equation}
Thanks to the Duhamel formula, we have
\begin{equation}
\D_l \ujt = e^{t\G} \D_l \D_j u_0 +\int_0^t e^{(t-\tau)\G}\D_l \bigg(\Fjt+\Gjt+\Rjt+\SjtL+\SjtNL\bigg)(\tau) d\tau,
\label{Duham}
\end{equation}
and then we can follow the very same lines as in \cite{TH1, Dlagrangien, CD, CHVP, Corder} provided that we are able to estimate $\|\D_l\SjtL\|_{L^p}$ and $\|\D_l\SjtNL\|_{L^p}$, which is the object of the following result:
\begin{prop}
\sl{Under the same assumptions, there exists two constants $C>0$ and $C_F>0$ such that, for all $p\in[1,\infty]$, if $t$ is so small that:
$$
e^{CV(t)}-1 \leq \frac{1}{2},
$$
then
$$
\|\D_l\SjtL\|_{L^p} +\|\D_l\SjtNL\|_{L^p} \leq C_{F}\max(\nu, \nu') 2^{3j-l} e^{CV(t)}(e^{CV(t)}-1) \|u_j(t)\|_{L^p},
$$
where $V(t)=\int_0^t \|\n v(\tau)\|_{L^\infty} d\tau$.
}
\label{estimcomm}
\end{prop}
\textbf{Proof :} The first estimate is easily dealt exactly as in the case of the classical Laplacian and we will skip details (we refer to \cite{TH1, FC3, CD}). Estimating the non-local commutator is the main difficulty of the present article and will require more attention. This is the object of the following section.

\subsection{Commutator estimates for $\SjtNL$}

First, as in \cite{TH3} we simplify the problem by rewriting $\SjtNL$ as follows:
\begin{multline}
(\La^2 \uj)\circ \psj-\La^2(\uj \circ \psj)=\big(\La(\La\uj)\big)\circ \psj-\La\big((\La\uj) \circ \psj\big)\\
+\La \bigg((\La \uj)\circ \psj-\La(\uj \circ \psj)\bigg),
\label{decpLa2}
\end{multline}
which allows us to reduce the study to the following quantity: for any function $f$,
$$
I_j=I_j(f)= \big(\La f\big) \circ \psj -\La \big(f \circ \psj\big).
$$
The first step in the study of $\SjtNL$ is then the following result, whose proof is given in the next section:
\begin{prop}
\sl{There exist two constants $C,C_F>0$ so that for all $p\in[1,\infty]$ and all function $f$, denoting $f_j =\D_j f$ and, as usual, $V(t)=\int_0^t \|\n v(\tau)\| d\tau$, if $t$ is so small that
$$
e^{2CV(t)}-1 \leq \frac{1}{2},
$$
then
$$
\|I_j(f_j)\|_{L^p} \leq C_F e^{CV(t)} (e^{CV(t)}-1) 2^j \|f_j\|_{L^p}.
$$
}
\label{propLa}
\end{prop}

\subsubsection{Proof of the commutator estimates for $\La$}

As in \cite{CHVP, Corder}, instead of estimating $I_j(x)$, it will be simpler to estimate $I_j(\psjm(x))$ as it considerably simplifies the integrals. Moreover, contrary to \cite{CHVP, Corder} we have $\div v=0$ so the Jacobian determinant is equal to $1$ and the previous quantities have the same $L^p$-norm. Most of the computations are valid for any $\cC^2$-function $f$, and we will precise when it will be needed that $f$ is spectrally localized in $2^j \cC$ (we will write $f_j$). Thanks to \eqref{exprdiff1}, for all $x\in \R^3$,
\begin{multline}
I_j(\psjm(x))= (\La f)(x) -\La \big(f \circ \psj\big)(\psjm(x))\\
=\lim_{\ee \rightarrow 0} \left( \int_{|y|\geq \ee}  K(y)\big(f(x-y)-f(x)\big) dy -\int_{|y|\geq \ee}  K(y)\bigg(\big(f\circ \psj\big)(\psjm(x)-y)-f(x)\bigg) dy \right).
\end{multline}
Performing in the second integral the following change of variable (see \cite{CHVP})
$$
x-z=\psj\big(\psjm(x)-y\big) \Leftrightarrow y=\psjm(x)-\psjm(x-z),
$$
we obtain that
\begin{equation}
I(\psjm(x))= \lim_{\ee \rightarrow 0} g_\ee(x),
\label{geps}
\end{equation}
where
\begin{multline}
g_\ee(x) =\int_{|y|\geq \ee}  K(y)\big(f(x-y)-f(x)\big) dy\\
-\int_{|m_x(-y)|\geq \ee}  K\big(m_x(-y)\big) \big(f(x-y)-f(x)\big) dy.
\end{multline}
with the following notations: for all $x,y\in \R^3$,
\begin{equation}
\begin{cases}
\vspace{1mm}
m_x(y)=\psjm(x)-\psjm(x+y),\\
\displaystyle{Y_{\pm}=\frac{|m_x(\pm y)|}{|y|} \quad \mbox{and}\quad Y_{\pm}^F=\frac{|m_x(\pm y)|_\frac{1}{F}}{|y|_\frac{1}{F}} \quad \mbox{if } y\neq 0.}
\end{cases}
\label{notmxY}
\end{equation}
In what follows we will make an extensive use of estimates involving these diffeomorphisms and for simplicity, we put in the appendix all these properties, that are much more precise than in \cite{CHVP, Corder}.\\

The main ingredient in the proof of Proposition \ref{propLa} is to rewrite $I_j(\psjm(x))$ in a more handy way, which is the object of the following result:
\begin{prop}
\sl{Under the same assumptions, for all $x\in \R^3$,
\begin{multline}
I_j(\psjm(x))=\frac{1}{2} \bigg( \int_{\R^3} \big(K(y)-K(m_x(-y)\big) \big(f(x-y)+f(x+y)-2f(x)\big)dy\\
+\int_{\R^3} \big(K(m_x(-y))-K(m_x(y))\big) \big(f(x+y)-f(x)\big) dy \bigg).
\end{multline}
\label{propIrewrite}
}
\end{prop}
\textbf{Proof :} For all $\ee>0$, and $x\in \R^3$,
\begin{multline}
g_\ee(x)= \int_{\R^3}  \bigg(K(y){\bf{1}}_{\{|y|\geq \ee\}} -K\big(m_x(-y)\big) {\bf{1}}_{\{|m_x(-y)|\geq \ee\}}\bigg) \big(f(x-y)-f(x)\big) dy\\
=\int_{\R^3}  \bigg(K(y){\bf{1}}_{|y|\geq \ee} \big({\bf{1}}_{\{|m_x(-y)|\geq \ee\}} +{\bf{1}}_{\{|m_x(-y)|< \ee\}}\big)\\
-K\big(m_x(-y)\big) {\bf{1}}_{\{|m_x(-y)|\geq \ee\}} \big({\bf{1}}_{\{|y|\geq \ee\}} +{\bf{1}}_{\{|y|< \ee\}}\big)\bigg) \big(f(x-y)-f(x)\big) dy.
\end{multline}
so that we can rewrite $g_\ee$ into:
$$
g_\ee(x) =g_\ee^1(x) +g_\ee^2(x) -g_\ee^3(x),
$$
with
\begin{equation}
\begin{cases}
\vspace{1mm}
\displaystyle{g_\ee^1(x) =\int_{\R^3} {\bf{1}}_{\{|y|\geq \ee\}} {\bf{1}}_{\{|m_x(-y)|\geq \ee\}} \bigg(K(y)-K\big(m_x(-y)\big)\bigg) \big(f(x-y)-f(x)\big) dy,}\\
\vspace{1mm}
\displaystyle{g_\ee^2(x) =\int_{\R^3} {\bf{1}}_{\{|y|\geq \ee\}} {\bf{1}}_{\{|m_x(-y)|< \ee\}} K(y) \big(f(x-y)-f(x)\big) dy,}\\
\displaystyle{g_\ee^3(x) =\int_{\R^3} {\bf{1}}_{\{|y|< \ee\}} {\bf{1}}_{\{|m_x(-y)|\geq \ee\}} K(m_x(-y)) \big(f(x-y)-f(x)\big) dy,}
\end{cases}
\label{g123}
\end{equation}
We will now prove the following result:
\begin{lem}
\sl{Under the previous assumptions,
$$
\sup_{x\in\R^3} \big(|g_\ee^2(x)| +|g_\ee^3(x)|\big) \underset{\ee \rightarrow 0}{\longrightarrow} 0.
$$
\label{lemg23}
}
\end{lem}
\textbf{Proof: } We will only concentrate on $g_\ee^2$ because, up to applying the diffeomorphism $m_x$, $g_\ee^3$ can be dealt exactly the same way. First we emphasize that any direct estimate will fail as it will only give that $|g_\ee^2(x)|$ is bounded as $\ee$ goes to zero. Indeed, thanks to Proposition \ref{propKymx} from the appendix, we have:
\begin{multline}
|g_\ee^2(x)| \leq C_F \int_{\ee \leq |y|< e^{CV(t)}\ee} \frac{1}{|y|^3} e^{6CV(t)}(e^{2CV(t)}-1) \|\n f\|_{L^\infty} dy\\
\leq C_F e^{6CV(t)}(e^{2CV(t)}-1)(e^{3CV(t)}-1)\|\n f\|_{L^\infty}.
\end{multline}
In order to prove the desired result we need, as in \eqref{exprdiff2}, to perform the change of variable $y\mapsto -y$ in the second term of the obvious identity $g_\ee^2(x)=\frac{1}{2}(g_\ee^2(x)+g_\ee^2(x))$ , which allows us to write that for all $x\in \R^3$:
$$
g_\ee^2(x) =\frac{1}{2} \bigg(\int_{A_\ee^-} K(y) \big(f(x-y)-f(x)\big) dy +\int_{A_\ee^+} K(y) \big(f(x+y)-f(x)\big) dy\bigg)
$$
where we have introduced the following sets (of course depending on $x$):
\begin{equation}
A_\ee^\pm= \{y\in \R^3/|y|\geq \ee \quad \mbox{and} \quad |m_x(\pm y)|< \ee\}.
\end{equation}
Obviously, we can state that:
\begin{lem}
\sl{The sets $A_\ee^\pm$ satisfy:
$$
\forall y\in \R^3, \quad y\in A_\ee^{-} \Leftrightarrow -y\in A_\ee^{+},
$$
and for $\eta\in\{-1,1\}$:
$$
A_\ee^\eta \subset \{y\in \R^3/\ee \leq |y|< \ee e^{CV(t)} \quad \mbox{and} \quad \ee e^{-CV(t)} \leq |m_x(\eta y)| <\ee\}.
$$
and
$$
\mbox{vol}(A_\ee^\eta) \leq \frac{4\pi}{3}\left(e^{3CV(t)}-1\right) \ee^3
$$
}
\label{lemA+}
\end{lem}
\textbf{Proof: } thanks to the first estimate from Proposition \ref{propmx}, we immediately obtain the second point. $\square$
\\

Taking these sets into consideration, we can rewrite $g_\ee^2$:
\begin{equation}
g_\ee^2(x) =\frac{1}{2} \left(I_\ee(x)+II_\ee(x)+III_\ee(x)\right),
\label{gee}
\end{equation}
with
$$
\begin{cases}
\vspace{1mm}
\displaystyle{I_\ee(x)= \int_{A_\ee^- \cap A_\ee^+} K(y) \big(f(x-y)+f(x+y)-2f(x)\big) dy,}\\
\vspace{1mm}
\displaystyle{II_\ee(x)= \int_{A_\ee^-\setminus A_\ee^+} K(y) \big(f(x-y)-f(x)\big) dy,}\\
\displaystyle{III_\ee(x)= \int_{A_\ee^+\setminus A_\ee^-} K(y) \big(f(x+y)-f(x)\big) dy.}
\end{cases}
$$
Lemma \ref{lemA+} immediately implies that $II_\ee=III_\ee$. We recall that we want to prove that $g_\ee^2$ goes to zero (see Lemma \ref{lemg23}). Let us begin with $I_\ee$: using twice the mean-value Theorem and Lemma \ref{lemA+} implies that for all $x\in \R^3$,
\begin{multline}
|I_\ee(x)| \leq \int_{\ee \leq |y| <\ee e^{CV(t)}} \frac{C_F}{|y|^4} |y|^2\|\n^2 f\|_{L^\infty} dy\\
\leq \frac{C_F}{\ee^2} \|\n^2f\|_{L^\infty} \mbox{vol} \left\{y\in \R^3, \quad \ee \leq |y| <\ee e^{CV(t)} \right\}\\
\leq C_F \|\n^2f\|_{L^\infty} \frac{4\pi}{3}\left(e^{3CV(t)}-1\right) \ee
\end{multline}
so that we obtain:
\begin{equation}
\sup_{x\in\R^3} |I_\ee(x)| \underset{\ee \rightarrow 0}{\longrightarrow} 0.
\label{Iee}
\end{equation}
As before, applying this rough argument to $g_\ee^2$ only provides that it is bounded whereas we need to prove it goes to zero. For $II_\ee$ (and $III_\ee$) as we only have a finite difference of order one, we need to be much more precise when estimating the volume of $A_\ee^-\setminus A_\ee^+$:
$$
|II_\ee(x)| \leq \int_{A_\ee^-\setminus A_\ee^+} \frac{C_F}{|y|^4} |y| \|\n f\|_{L^\infty} dy \leq \frac{C_F}{\ee^3} \|\n f\|_{L^\infty} \mbox{vol} \left(A_\ee^-\setminus A_\ee^+ \right).
$$
In fact, as each of the sets $A_\ee^\pm$ is already very close to the domain located between the spheres centered at zero and with radii $\ee$ and $\ee e^{CV(t)}$, their symmetric difference has an even smaller volume:
\begin{lem}
\sl{Under the previous assumptions and notations, there exists a constant $C_{F,j}$ (also depending on $t$ but bounded) such that:
$$
\mbox{vol}(A_\ee^-\setminus A_\ee^+) \leq C_{F,j} \ee^4
$$
}
\label{lemvoldiffA}
\end{lem}
\textbf{Proof :} Let us write that
$$
A_\ee^-\setminus A_\ee^+ \subset\left\{y\in \R^3, \quad \ee \leq |y| <\ee e^{CV(t)},\quad \ee e^{-CV(t)} \leq |m_x(-y)| <\ee ,\mbox{ and} \quad |m_x(y)| \geq\ee \right\}.
$$
Thanks to the last point of Proposition \ref{propmx}, if $y\in A_\ee^-\setminus A_\ee^+$ then:
\begin{multline}
\bigg||m_x(-y)|-|m_x(y)|\bigg| \leq  e^{2C V(t)}(e^{2C V(t)}-1) \min(1,2^j|y|)|y|\\
\leq e^{4C V(t)}(e^{2C V(t)}-1) 2^j \ee^2,
\end{multline}
so that, using this with the fact that $|m_x(-y)| <\ee$, we obtain
$$
\ee \leq |m_x(y)| \leq |m_x(-y)| +\bigg||m_x(-y)|-|m_x(y)|\bigg| <\ee\left(1+e^{4C V(t)}(e^{2C V(t)}-1) 2^j \ee\right).
$$ 
That is
$$
A_\ee^-\setminus A_\ee^+ \subset\left\{y\in \R^3, \quad  \ee \leq |m_x(y)| <\ee\left(1+e^{4C V(t)}(e^{2C V(t)}-1) 2^j \ee\right)\right\}.
$$
Then, as $y\mapsto m_x(y)$ is a volume preserving diffeomorphism, roughly denoting $C_j=e^{4C V(t)}(e^{2C V(t)}-1) 2^j$, we get:
\begin{multline}
\mbox{vol} \left(A_\ee^-\setminus A_\ee^+ \right) \leq \mbox{vol} \left\{ y\in \R^3, \quad  \ee \leq |y| <\ee\left(1+ \ee C_j\right)\right\}\\
=\frac{4\pi}{3} \ee^3\left((1+\ee C_j)^3-1\right) =\frac{4\pi}{3} \ee^4\left(3 C_j+ 3\ee C_j^2 + \ee^2 C_j^3\right),
\end{multline}
which concludes the proof of the lemma. $\square$
\begin{rem}
\sl{We do not need to worry about the dependency in time, as in what follows we will deal with small $t$ such that $e^{CV(t)}-1 \leq 1/2$.
}
\end{rem}
This allows us to conclude that there exists a constant $C_{F,j}>0$ such that:
\begin{equation}
\sup_{x\in\R^3} |II_\ee(x)| \leq \ee C_{F,j}\|\n f\|_{L^\infty} \underset{\ee \rightarrow 0}{\longrightarrow} 0.
\label{IIee}
\end{equation}
Gathering \eqref{gee}, \eqref{Iee} and \eqref{IIee} ends the proof of Lemma \ref{lemg23}. $\square$

This immediately leads us to: 
\begin{multline}
I_j(\psjm(x))=\lim_{\ee \rightarrow 0} g_\ee(x) =\lim_{\ee \rightarrow 0} g_\ee^1(x)\\
=\lim_{\ee \rightarrow 0} \int_{\R^3} {\bf{1}}_{\{|y|\geq \ee\}} {\bf{1}}_{\{|m_x(-y)|\geq \ee\}} \bigg(K(y)-K(m_x(-y))\bigg) \big(f(x-y)-f(x)\big) dy.
\end{multline}
We emphasize that, for the same reason as before, we still cannot let $\ee$ go to zero in this integral: we obtain thanks to Proposition \ref{propKymx}, that near zero, the integrated function is bounded by:
$$
\frac{1}{|y|^3}(e^{2CV(t)}-1)\|\n f\|_{L^\infty},
$$
which is not integrable at $0$. Our only option is to use once more the same desingularization argument as before: $g_\ee^1(x)= \frac{1}{2} (g_\ee^1(x)+g_\ee^1(x))$ and perform the change of variable $y\mapsto -y$ in the second term, in order to take advantage of finite differences of order $2$. Doing as what lead to \eqref{g123} we can write that:
$$
g_\ee^1(x)= IV_\ee(x) +V_\ee(x) +VI_\ee(x),
$$
where
$$
IV_\ee(x) =\int_{B_\ee} \bigg(\big(K(y)-K(m_x(-y)\big) \big(f(x-y)-f(x)\big) +\big(K(y)-K(m_x(y)\big) \big(f(x+y)-f(x)\big) \bigg)dy,
$$
with
$$
B_\ee=\{y\in \R^3, \quad |y|\geq \ee,\quad |m_x(y)| \geq \ee,\quad |m_x(-y)| \geq \ee \},
$$
and
\begin{equation}
\begin{cases}
\vspace{1mm}
\displaystyle{V_\ee(x) =\int_{A_\ee^+ \setminus A_\ee^-} \bigg(K(y)-K(m_x(-y)\bigg) \big(f(x-y)-f(x)\big) dy,}\\
\vspace{1mm}
\displaystyle{VI_\ee(x) =\int_{A_\ee^- \setminus A_\ee^+} \bigg(K(y)-K(m_x(y)\bigg) \big(f(x+y)-f(x)\big) dy}
\end{cases}
\end{equation}
As before, thanks to the change of variable $y\mapsto -y$, we have $V_\ee(x)=VI_\ee(x)$, and thanks to Proposition \ref{propKymx} (see Appendix) and Lemma \ref{lemvoldiffA},
\begin{multline}
\sup_{x\in\R^3} |V_\ee(x)| \leq \int _{A_\ee^+ \setminus A_\ee^-} \frac{C_F}{|y|^3} e^{6CV(t)} \big(e^{6CV(t)} -1\big) \|\n f\|_{L^\infty} dy\\
\leq \frac{C_F}{\ee^3} e^{6CV(t)} \big(e^{6CV(t)} -1\big) \|\n f\|_{L^\infty} \mbox{vol}\big( A_\ee^+ \setminus A_\ee^-\big)\\
\leq C_{F,j} e^{6CV(t)} \big(e^{6CV(t)} -1\big) \|\n f\|_{L^\infty} \ee \underset{\ee \rightarrow 0}{\longrightarrow} 0.
\end{multline}
From this we finally deduce that:
\begin{multline}
I_j(\psjm(x))=\frac{1}{2} \lim_{\ee \rightarrow 0} IV_\ee(x) =\frac{1}{2} \lim_{\ee \rightarrow 0} \bigg[ \int_{B_\ee} \bigg(K(y)-K(m_x(-y)\bigg) \big(f(x-y)+f(x+y)-2f(x)\big)dy\\
+\int_{B_\ee} \bigg(K(m_x(-y))-K(m_x(y))\bigg) \big(f(x+y)-f(x)\big) dy \bigg] =\frac{1}{2} \left[\A(x) +\B(x)\right],
\end{multline}
where
\begin{equation}
\begin{cases}
\vspace{1mm}
\A(x)= \displaystyle{\int_{\R^3} \bigg(K(y)-K(m_x(-y)\bigg) \big(f(x-y)+f(x+y)-2f(x)\big)dy,}\\
\B(x)= \displaystyle{\int_{\R^3} \bigg(K(m_x(-y))-K(m_x(y))\bigg) \big(f(x+y)-f(x)\big) dy \bigg).}
\end{cases}
\end{equation}
This concludes the proof of Proposition \ref{propIrewrite}. $\blacksquare$
\begin{rem}
\sl{We emphasize that, at last, we could pass to the limit $\ee \rightarrow 0$ in the integrals as the functions are now integrable on $\R^3$.}
\end{rem}
\begin{rem}
\sl{Note that this result was completely trivial if we could pass to the limit in \eqref{exprdiff1} and \eqref{geps}: indeed if the functions were integrable we only had to perform the usual change of variable $y\mapsto -y$.}
\end{rem}
To conclude the proof of Proposition \ref{propLa}, we will study $\A$ and $\B$: the first term is very easy to estimate. Thanks to Proposition \ref{propKymx}, and Theorem \ref{thdiff2}, if $e^{CV(t)}-1\leq 1/2$,
\begin{multline}
\|\A\|_{L^p} \leq \int_{\R^3} \frac{C_F}{|y^4|} e^{6CV(t)} \big(e^{2CV(t)}-1\big) \|f_j(.-y)+f_j(.+y)-2f_j(.)\|_{L_x^p} dy\\
\leq C_F e^{6CV(t)} \big(e^{2CV(t)}-1\big) \|f_j\|_{\dot{B}_{p,1}^1} \leq  C_F e^{6CV(t)} \big(e^{2CV(t)}-1\big) 2^j \|f_j\|_{L^p}.
\end{multline}
\begin{rem}
\sl{Note that Theorem \ref{thdiff2} provides estimates for homogeneous Besov spaces, but as we use these results with $f_j$, which is spectrally localized, this entails the result.}
\end{rem}
In order to estimate correctly $\B$ we need a careful estimate of $K(m_x(-y))-K(m_x(y))$. We refer to the appendix for the proof but we emphasize that, as in \cite{CHVP, Corder} this term provides the required power of $|y|$ allowing the function to be integrable at $0$. More precisely, thanks to Proposition \ref{propKymx},
\begin{multline}
\|\B\|_{L^p} \leq \int_{\R^3} \frac{C_F}{|y|^4} e^{6CV(t)} \big(e^{2CV(t)}-1\big) \min(1, 2^j|y|) \|f_j(.+y)-f_j(.)\|_{L_x^p} dy\\
\leq C_F e^{6CV(t)} \big(e^{2CV(t)}-1\big) \int_{\R^3} \frac{\min(1, 2^j|y|)}{|y|^4} \|f_j(.+y)-f_j(.)\|_{L_x^p} dy.
\end{multline}
The rest of the proof is classical: we refer to \cite{Dbook} (section $2.4$, se also \cite{CHVP}) for the following estimate:
$$
\|f_j(.+y)-f_j(.)\|_{L^p} \leq C \min(1, 2^j|y|) \|f_j\|_{L^p},
$$
so that, plugging this in the previous integral gives:
\begin{multline}
\|\B\|_{L^p} \leq C_F e^{6CV(t)} \big(e^{2CV(t)}-1\big) \|f_j\|_{L^p} \int_{\R^3} \frac{\min(1, 2^j|y|)^2}{|y|^4} dy\\
\leq C_F e^{6CV(t)} \big(e^{2CV(t)}-1\big) \|f_j\|_{L^p} \left( \int_{|y|\leq 2^{-j}} \frac{2^{2j}}{|y|^2} dy + \int_{|y|\geq 2^{-j}} \frac{1}{|y|^4} dy \right)\\
\leq  C_F e^{6CV(t)} \big(e^{2CV(t)}-1\big) 2^j \|f_j\|_{L^p}.
\end{multline}
This concludes the proof of Proposition \ref{propLa}. $\blacksquare$

\begin{rem}
\sl{In the case $p=\infty$ it is possible to prove the result in a slightly simpler way using the formulation with finite differences of order 2.}
\end{rem}

\subsubsection{Commutation for $\SjtNL$, end of the proof of Proposition \ref{estimcomm}}

We studied the commutation with $\La$ but we recall that our aim is to study the commutation with $\La^2$. This part of the proof is classical so we will skip details and only point out what changes. Thanks to \ref{decpLa2}, we can write that:
\begin{multline}
\|\D_l \bigg(\La^2 \uj)\circ \psj-\La^2(\uj \circ \psj\bigg)\|_{L^p} \leq \|\D_l\bigg(\big(\La(\La\uj)\big)\circ \psj-\La\big((\La\uj) \circ \psj\big)\bigg)\|_{L^p}\\
+\|\D_l \La \bigg((\La \uj)\circ \psj-\La(\uj \circ \psj)\bigg)\|_{L^p} = \|\D_l \mathbbm{1}\|_{L^p}+\|\D_l \mathbbm{2}\|_{L^p},
\label{decp12}
\end{multline}
Let us begin with $\mathbbm{1}$. Of course, if we roughly estimate this term, we end up with terms that are not summable (as we have to perform a summation over $l\geq j-N_0$). To bypass this problem, we simply use the same idea as Vishik (see \cite{Vishik, TH1, CD, CHVP, Corder}) and write (we will only use it in the case $l\geq 0$):
$$
\|\D_l \mathbbm{1}\|_{L^p} \leq 2^{-l} \|\D_l \n \mathbbm{1}\|_{L^p} \leq 2^{-l} \|\n \mathbbm{1}\|_{L^p},
$$
with, denoting $g_j= \La\uj$,
$$
\n \mathbbm{1}= \big(\La\n g_j\big)\circ \psj\cdot D\psj-\La\big(\n g_j \circ \psj \cdot D\psj\big) =\mathbbm{1}_A +\mathbbm{1}_B,
$$
where (we recall that $I_3$ denotes the unit matrix of size $3$ and $D \psj$ the jacobian matrix of $\psj$):
\begin{equation}
\begin{cases}
\vspace{1mm}
\mathbbm{1}_A =\big(\La\n g_j\big)\circ \psj\cdot (D\psj -I_3)-\La\big(\n g_j \circ \psj \cdot (D\psj-I_3)\big),\\
\mathbbm{1}_B =\big(\La\n g_j\big)\circ \psj-\La\big(\n g_j \circ \psj\big).
\end{cases}
\end{equation}
Thanks to Proposition \ref{propLa}, we obtain that (if $t$ is small enough):
\begin{equation}
\|\mathbbm{1}_B\|_{L^p} \leq C_F e^{CV(t)} (e^{CV(t)}-1) 2^j \|\n g_j\|_{L^p} \leq C_F e^{CV(t)} (e^{CV(t)}-1) 2^{3j} \|u_j\|_{L^p}.
\label{est1a}
\end{equation}
For the other term, introducing $k_j =D\psj -I_3$, we can write that $\mathbbm{1}_A =\mathbbm{1}_{A1} +\mathbbm{1}_{A2}$, with
\begin{equation}
\begin{cases}
\mathbbm{1}_{A1} =\bigg(\big(\La\n g_j\big)\circ \psj-\La\big(\n g_j \circ \psj \big)\bigg)\cdot k_j,\\
\mathbbm{1}_{A2} =\La \big(\n g_j \circ \psj \big)\cdot k_j-\La\big(\n g_j \circ \psj \cdot k_j\big).
\end{cases}
\end{equation}
The first term is easily dealt thanks to Propositions \ref{propLa} and \ref{p:flow}:
\begin{multline}
\|\mathbbm{1}_{A1} \|_{L^p} \leq \|\big(\La\n g_j\big)\circ \psj-\La\big(\n g_j \circ \psj \big)\|_{L^p}\cdot \|k_j\|_{L^\infty}\\
\leq C_F e^{CV(t)} (e^{CV(t)}-1)^2 2^j \|\n g_j\|_{L^p} \leq C_F e^{2 CV(t)} (e^{CV(t)}-1) 2^{3j} \|u_j\|_{L^p}.
\label{estim1a1}
\end{multline}
Estimating the second term will require more efforts and we begin with the following property of the non-local operator $\La$:
\begin{prop}
\sl{For any smooth functions $f,g$ we can write:
$$
\La (fg)= f\La g +g\La f +M(f,g),
$$
where the bilinear operator $M$ is defined by:
$$
M(f,g)(x) =\int_{\R^3} K(y) \big(f(x-y)-f(x)\big) \big(g(x-y)-g(x)\big) dy.
$$
Moreover there exists a constant $C_F$ such that for all $f,g$:
\begin{equation}
\|M(f,g)\|_{L^p} \leq C_F \sqrt{\|f\|_{L^p} \|\n f\|_{L^p} \|g\|_{L^\infty} \|\n g\|_{L^\infty}}.
\end{equation}
}
\label{propM}
\end{prop}
\textbf{Proof :} For the first part of the result, we simply use \eqref{exprdiff1}: for all $x\in\R^3$,
\begin{multline}
\big(\La (fg)-f\La g -g\La f \big)(x) =\lim_{\ee \rightarrow 0} \int_{|y|\geq \ee}  K(y)\bigg(\big((fg)(x-y)-(fg)(x)\big) -f(x)\big(g(x-y)-g(x)\big)\\
-g(x)\big(f(x-y)-f(x)\big) \bigg) dy\\
=\lim_{\ee \rightarrow 0} \int_{|y|\geq \ee}  K(y)\big(f(x-y)-f(x)\big) \big(g(x-y)-g(x)\big) dy.
\end{multline}
The integral converges on $\R^3$ (no problem at zero or at infinity) so we obtain the conclusion thanks to the Lebesgue theorem. For the estimate, we use a classical threshold argument: for a certain $R$ (to be fixed later), we can write that
\begin{multline}
\|M(f,g)\|_{L^p} \leq C_F \bigg( \int_{|y|\leq R}\frac{1}{|y|^4} \|\big(f(.-y)-f(.)\big) \big(g(.-y)-g(.)\big)\|_{L_x^p} dy\\
+\int_{|y|\geq R}\frac{1}{|y|^4} \|\big(f(.-y)-f(.)\big) \big(g(.-y)-g(.)\big)\|_{L_x^p} dy \bigg)\\
\leq C_F \bigg( \int_{|y|\leq R}\frac{1}{|y|^2} \|\n f\|_{L^p} \|\n g\|_{L^\infty} dy +\int_{|y|\geq R}\frac{1}{|y|^4} 4\|f\|_{L^p} \|g\|_{L^\infty} dy \bigg)\\
\leq C_F \left( R \|\n f\|_{L^p} \|\n g\|_{L^\infty} +\frac{1}{R}\|f\|_{L^p} \|g\|_{L^\infty} \right).
\end{multline}
Then, classically, this summation of two terms with a constant product is minimal when they are equal, that is when we choose:
$$
R^2=\frac{\|f\|_{L^p} |g\|_{L^\infty}}{\|\n f\|_{L^p} \|\n g\|_{L^\infty}}. \blacksquare
$$
We can go back to the estimation of $\mathbbm{1}_{A2}$: applying the previous estimate with $f= \n g_j \circ \psj$ and $g= k_j$, we obtain
\begin{multline}
\|\mathbbm{1}_{A2}\|_{L^p}= \|-M(\n g_j \circ \psj, k_j)+ (\n g_j \circ \psj) \cdot \La k_j\|_{L^p}\\
\leq C_F \bigg[ \bigg( \|\n g_j \circ \psj\|_{L^p} \|\n^2 g_j \circ \psj \cdot D\psj\|_{L^p} \|D\psj-I_3\|_{L^\infty} \|D^2 \psj\|_{L^\infty}\bigg)^{\frac{1}{2}}\\
+\|\n g_j \circ \psj\|_{L^p} \|\La k_j\|_{L^\infty} \bigg]\\
\end{multline}
Thanks to Proposition \ref{p:flow} and the fact that $u_j$ is frequency localized, we deduce
$$
\|\mathbbm{1}_{A2}\|_{L^p} \leq C_F \Bigg(2^{6j} \|u_j\|_{L^p}^2 e^{CV(t)} (e^{CV(t)}-1)^2 \bigg)^{\frac{1}{2}} +C_F 2^{2j} \|u_j\|_{L^p} \|\La k_j\|_{L^\infty}.
$$
We have to be particularly careful with the last term, as $k_j$ is not frequency localized. A way to deal with it is to write (see \eqref{estimLaBesov}):
$$
\|\La k_j\|_{L^\infty} =\|\La (D\psj-I_3)\|_{L^\infty} \leq C_F \|D\psj-I_3\|_{\dot{B}_{p,1}^1}.
$$
Then, we use the following estimate:
\begin{lem}
\sl{There exists a constant $C>0$ such that for all smooth function $u$, we have:
$$
\|u\|_{\dot{B}_{r,1}^1} \leq C \sqrt{\|u\|_{L^r} \|\n^2 u\|_{L^r}}.
$$
}
\end{lem}
\textbf{Proof :} it is the same idea as in Proposition \ref{propM}: for a certain $\aa$ to be fixed later, we have
\begin{multline}
\|u\|_{\dot{B}_{r,1}^1} \leq \Sum_{j\leq \aa} 2^j \|\ddj u\|_{L^r} + \Sum_{j\geq \aa+1} 2^{-j} \|\ddj \n^2 u\|_{L^r} \\
\leq \big(\Sum_{j\leq \aa} 2^j \big) \|u\|_{L^r} + \big(\Sum_{j\geq \aa+1} 2^{-j}\big) \|\n^2 u\|_{L^r} \\
\leq C \big(2^\aa \|u\|_{L^r} + 2^{-\aa} \|\n^2 u\|_{L^r} \big),
\end{multline}
and choosing $\aa$ the closest integer so that $2^{2\aa} \sim\|u\|_{L^r} /\|\n^2 u\|_{L^r}$ ends the proof. $\square$
\\
From that we deduce (using Proposition \ref{p:flow}):
$$
\|D\psj-I_3\|_{\dot{B}_{p,1}^1} \leq C \bigg( \|D\psj-I_3\|_{L^\infty} \|D^3 \psj\|_{L^\infty} \bigg)^{\frac{1}{2}} \leq C (e^{CV(t)}-1) 2^j,
$$
so that we can conclude:
\begin{equation}
\|\mathbbm{1}_{A2}\|_{L^p} \leq C_F 2^{3j} e^{CV(t)} (e^{CV(t)}-1) \|u_j\|_{L^p}.
\label{estim1a2}
\end{equation}
Gathering \eqref{estim1a1} and \eqref{estim1a2}, we obtain that:
\begin{equation}
\|\D_l \mathbbm{1}\|_{L^p} \leq C_F 2^{3j-l} e^{CV(t)} (e^{CV(t)}-1) \|u_j\|_{L^p}.
\label{estim11}
\end{equation}
\begin{rem}
\sl{
As a consequence of the previous lemma, we have:
\begin{multline}
\|\La(fg)-(\La f)g\|_{L^p} =\|f\La g+M(f,g)\|_{L^p}\\
\leq C_F \bigg( \|f\|_{L^p} \big(\|g\|_{L^\infty} \|\n^2 g\|_{L^\infty}\big)^\frac{1}{2}+\big(\|f\|_{L^p} \|\n f\|_{L^p} \|g\|_{L^\infty} \|\n g\|_{L^\infty}\big)^\frac{1}{2} \bigg),\\
\leq C_F \|f\|_{L^p}^\frac{1}{2} \|g\|_{L^\infty}^\frac{1}{2} \bigg( \|f\|_{L^p}^\frac{1}{2} \|\n^2 g\|_{L^\infty}^\frac{1}{2} +\|\n f\|_{L^p}^\frac{1}{2} \|\n g\|_{L^\infty}^\frac{1}{2} \bigg).
\label{commult}
\end{multline}
}
\end{rem}
We now turn to the second term $\mathbbm{2}$ from \eqref{decp12}. Exactly as in \cite{CHVP, Corder}, as we aim to obtain a result summable for $l\geq j-N_0$, we have to use twice the argument of Vishik on $\mathbbm{2}$ because of the additionnal $\La$:
\begin{multline}
\|\mathbbm{2}\|_{L^p}  =\|\D_l \La \bigg((\La \uj)\circ \psj-\La(\uj \circ \psj)\bigg)\|_{L^p} \leq 2^{-2l} \|\n^2\La \D_l\bigg((\La \uj)\circ \psj-\La(\uj \circ \psj)\bigg)\|_{L^p}\\
\leq 2^{-l} \|\n^2 \bigg((\La \uj)\circ \psj-\La(\uj \circ \psj)\bigg)\|_{L^p},
\end{multline}
and the rest of the proof is similar, we first decompose:
$$
\n^2 \bigg((\La \uj)\circ \psj-\La(\uj \circ \psj)\bigg) =\mathbbm{2}_A +\mathbbm{2}_B +\mathbbm{2}_C +\mathbbm{2}_D, 
$$
where
$$
\begin{cases}
\mathbbm{2}_A =(\La \n^2\uj)\circ \psj\cdot D\psj(D\psj-I_3)-\La\bigg(\n^2\uj \circ \psj \cdot D\psj(D\psj-I_3\bigg),\\
\mathbbm{2}_B =(\La \n\uj)\circ \psj \cdot D^2\psj-\La\bigg(\n\uj \circ \psj \cdot D^2\psj\bigg),\\
\mathbbm{2}_C =(\La \n^2\uj)\circ \psj\cdot (D\psj-I_3)-\La\bigg(\n^2\uj \circ \psj \cdot (D\psj-I_3)\bigg),\\
\mathbbm{2}_D =(\La \n^2\uj)\circ \psj-\La(\n^2\uj \circ \psj).
\end{cases}
$$
Using the very same arguments as for $\mathbbm{1}_{A2}$ (we skip details) we obtain that:
\begin{equation}
\|\D_l \mathbbm{2}\|_{L^p} \leq C_F 2^{3j-l} e^{CV(t)} (e^{CV(t)}-1) \|u_j\|_{L^p}.
\label{estim12}
\end{equation}
Gathering \eqref{estim11} and \eqref{estim12} ends the proof of Proposition \ref{estimcomm}. $\blacksquare$

\section{Proof of Theorem \ref{thC2}}

We can now return to the proof of the a priori estimates. Let us begin with Theorem \ref{thC2}. Without any loss of generality, we will prove it for $T_1=0$ and $T_2=T$. We recall that, for the high frequencies, thanks to \eqref{Duham} and Proposition \ref{semig2}, there exist two constants $C,c>0$ such that:
\begin{multline}
\|\D_l \ujt\|_{L^p} \leq C\bigg[e^{-c \no t 2^{2l}} \|\D_l \D_j u_0\|_{L^p} \\
+\int_0^t e^{-c \no (t-\tau) 2^{2l}} \bigg(\|\D_l \Fjt\|_{L^p}+\|\D_l \Gjt\|_{L^p}+\|\D_l \Rjt\|_{L^p}+\|\D_l (\SjtL+\SjtNL)\|_{L^p} \bigg)(\tau) d\tau \bigg].
\label{estimLptilde}
\end{multline}
First, we refer to \cite{TH1} for the following result:
\begin{prop}
\sl{
With the notations from \eqref{Rj}, there exists a constant $C>0$ such that:
$$
\|R_j\|_{L^p} \leq C \|\n v\|_{L^\infty} \|u\|_{L^p}.
$$
}
\label{propRjT1}
\end{prop}
Using once again the method of Vishik, we deduce that (we refer to the previous section for details):
\begin{multline}
\|\D_l \Rjt\|_{L^p} \leq C 2^{-l} \|\n (\Rjt\circ \psj)\|_{L^p} \leq C 2^{j-l} \|R_j\|_{L^p} e^{CV(t)}\\
\leq C 2^{j-l} e^{CV(t)}\|\nabla v(t)\|_{L^\infty} \|u(t)\|_{L^p}.
\end{multline}
Thanks to \eqref{Linftpetit}, as $\|v\|_{L_t^\infty L^6}\leq C'$, and if $t\in [0,T]$ for $T>0$ so small that $CC'T^\frac{1}{4}\leq \frac{1}{2} \no^\frac{3}{4}$, then we can write that:
\begin{equation}
\|\D_l \Rjt\|_{L^p} \leq C 2^{j-l} e^{CV(t)}\|\nabla v(t)\|_{L^\infty} \left(\|u_0\|_{L^p}+\int_0^t \big(\|\Fe(\tau)\|_{L^p} +\|\Ge(\tau)\|_{L^p}\big)d\tau\right).
\label{estimRj}
\end{equation}
\begin{rem}
\sl{Remember that in the setting of Theorem \ref{thC2}, $\Ge=0$.}
\end{rem}
Next, we obtained in Proposition \ref{estimcomm} that if $t$ is so small that 
$$
e^{CV(t)}-1 \leq \frac{1}{2},
$$
then
$$
\|\D_l\SjtL\|_{L^p} +\|\D_l\SjtNL\|_{L^p} \leq C_F\max(\nu, \nu') 2^{3j-l} e^{CV(t)}(e^{CV(t)}-1) \|u_j(t)\|_{L^p}.
$$
Then the arguments are the very same as in \cite{TH1, FC3} so we will skip details and refer the reader to the proofs and lemmas therein. Taking the $L^r$-norm in time ($r\in[1,\infty]$), using \eqref{estimRj}, \eqref{estimcomm} and the Young estimates on \eqref{estimLptilde}, we obtain that for all $t\in[0,T]$:
\begin{multline}
\|\D_l \ujt\|_{L_t^r L^p} \leq C(\no r 2^{2l})^{-\frac{1}{r}} \bigg[\|\D_l \D_j u_0\|_{L^p} +\int_0^t \|\D_l \Fjt(\tau)\|_{L^p} d\tau\\
+2^{j-l} e^{CV(t)} V(t) \left(\|u_0\|_{L^p}+\int_0^t \|\Fe(\tau)\|_{L^p} d\tau\right)\bigg]\\
+(\no r 2^{2l})^{-1} C_F\max(\nu, \nu') 2^{3j-l} e^{CV(t)}(e^{CV(t)}-1) \|u_j(t)\|_{L_t^r L^p}.
\end{multline}
After using that $\|\D_l \Fjt\|_{L^p} \leq C\|\Fjt\|_{L^p}\leq C\|\Fe\|_{L^p}$ we sum over $l\geq j-N_0$ ($N_0$ will be fixed soon), and get that for all $j\geq N_0$ and $t\in[0,T]$:
\begin{multline}
(\no r 2^{2j})^\frac{1}{r} \Sum_{l\geq j-N_0} \|\D_l \ujt\|_{L_t^r L^p} \leq C \|\D_j u_0\|_{L^p} +C_r 2^{2N_0}\int_0^t \|\Fe(\tau)\|_{L^p} d\tau\\
+C 2^{3N_0} e^{CV(t)} (e^{CV(t)}-1) \left(\|u_0\|_{L^p}+\int_0^t \|\Fe(\tau)\|_{L^p} d\tau\right)\\
+C_F\frac{\max(\nu, \nu')}{\no} 2^{3N_0} e^{CV(t)}(e^{CV(t)}-1) (\no r 2^{2j})^\frac{1}{r} \|u_j\|_{L_t^r L^p}.
\label{estimHF}
\end{multline}
Next, we turn to the low frequencies ($l\leq j-N_0$) which we can estimate as follows (we refer to \cite{Dbook}, Lemma $2.6$ \cite{TH1}, section $2.4.1$ in the divergence-free case), there exists a constant $C>0$ such that :
\begin{equation}
\|S_{j-N_0}\ujt\|_{L_t^r L^p} \leq C 2^{-N_0} e^{CV(t)} \|u_j\|_{L_t^r L^p}.
\label{estimBFR}
\end{equation}
Then, plugging this into \eqref{decpjl}, and gathering with \eqref{estimHF}, we can write that for all $j\geq N_0$:
\begin{multline}
(\no r 2^{2j})^\frac{1}{r} \|u_j\|_{L_t^r L^p} \leq C \|\D_j u_0\|_{L^p} +C_r 2^{2N_0}\int_0^t \|\Fe(\tau)\|_{L^p} d\tau\\
+\bigg( C2^{-N_0} e^{CV(t)} + (e^{CV(t)}-1)C_F\frac{\max(\nu, \nu')}{\no} 2^{3N_0} e^{CV(t)}\bigg) (\no r 2^{2j})^\frac{1}{r} \|u_j\|_{L_t^r L^p}\\
+C 2^{3N_0} e^{CV(t)} (e^{CV(t)}-1) \left(\|u_0\|_{L^p}+\int_0^t \|\Fe(\tau)\|_{L^p} d\tau\right).\\
\end{multline}
This is here that we fix $N_0$: we need $T>0$ small enough and $N_0$ large enough so that:
$$
\begin{cases}
\vspace{1mm}
\displaystyle{CC'T^\frac{1}{4} \leq \frac{1}{2} \no^\frac{3}{4},}\\
\vspace{1mm}
\displaystyle{e^{CV(T)}-1 \leq \frac{1}{2},}\\
\displaystyle{C2^{-N_0} e^{CV(T)} + (e^{CV(T)}-1)C_F\frac{\max(\nu, \nu')}{\min(\nu, \nu')} 2^{3N_0} e^{CV(T)} \leq \frac{1}{2},}
\end{cases}
$$ 
which is satified if:
$$
\begin{cases}
\vspace{1mm}
\displaystyle{T \leq \frac{\no^3}{(2CC')^4},}\\
\vspace{1mm}
\displaystyle{e^{CV(T)}\leq \frac{3}{2},}\\
\displaystyle{\frac{3}{2}C2^{-N_0} +(e^{CV(T)}-1)C_F\frac{\max(\nu, \nu')}{\min(\nu, \nu')} 2^{3N_0} \frac{3}{2} \leq \frac{1}{2},}
\end{cases}
$$ 
so that we first choose $N_0 \in\N$ large enough so that $2^{N_0} \geq 6C$ and then we take $T>0$ so small that:
\begin{equation}
\begin{cases}
\vspace{1mm}
\displaystyle{T \leq \frac{\no^3}{(2CC')^4},}\\
\vspace{1mm}
\displaystyle{e^{CV(T)}-1 \leq \frac{1}{6C_F\frac{\max(\nu, \nu')}{\min(\nu, \nu')} 2^{3N_0}}.}\\
\end{cases}
\label{condN0t}
\end{equation}
With these assumptions, we immediately obtain that there exists a constant $C_F$ such that for all $j\geq N_0$,
\begin{multline}
(\no r 2^{2j})^\frac{1}{r} \|u_j\|_{L_t^r L^p}\\
\leq C \|\D_j u_0\|_{L^p} +C_r \int_0^t \|\Fe(\tau)\|_{L^p} d\tau +C_F \left(\|u_0\|_{L^p}+\int_0^t \|\Fe(\tau)\|_{L^p} d\tau\right).
\end{multline}
For the low frequencies, thanks to Proposition \ref{estimLp} and condition \eqref{condN0t}, we simply write that for all $j\leq N_0$ (with the usual changes when $r=\infty$),
\begin{multline}
(\no r 2^{2j})^\frac{1}{r} \|u_j\|_{L_t^r L^p} \leq (\no r 2^{2N_0})^\frac{1}{r} \|u\|_{L_t^\infty L^p} t^\frac{1}{r}\\
\leq C(\no r)^\frac{1}{r} \left(\|u_0\|_{L^p}+\int_0^t \|\Fe(\tau)\|_{L^p} d\tau\right).
\end{multline}
Combining the last two estimates ends the proof of Theorem \ref{thC2}.

\section{Proof of Theorem \ref{thCs}}

The proof begins as for Theorem \ref{thC2}: we start from \eqref{estimLptilde}. There is no change for the non-local commutators $\D_l (\SjtL+\SjtNL)$, and for the remainder term $R_j$, instead of using Proposition \ref{propRjT1}, we will use the following result (we refer to \cite{TH1}, Lemma $2.4.1$ for the proof):
\begin{prop}
\sl{
With the same notations as in \eqref{Rj}, there exists a constant $C>0$ such that:
$$
\|R_j\|_{L^p} \leq C \|\n v\|_{L^\infty} \Sum_{k\geq -1}2^{-|k-j|}\|\D_k u\|_{L^p}.
$$
}
\label{propRjT2}
\end{prop}
Thanks to Lemma $2.6$ and $2.7$ from \cite{Dbook} (in the volume preserving case), we get:
$$
\begin{cases}
\vspace{1mm}
\|\D_l \Fjt(\tau)\|_{L^p} \leq C 2^{-|l-j|} e^{CV(\tau)} \|\Fj(\tau)\|_{L^p},\\
\displaystyle{\|\D_l \Rjt(\tau)\|_{L^p} \leq C 2^{-|l-j|} e^{CV(\tau)} \|\n v(\tau)\|_{L^\infty} \Sum_{k\geq -1}2^{-|k-j|}\|\D_k u(\tau)\|_{L^p}.}
\end{cases}
$$
From this we deduce that:
\begin{multline}
\|\D_l \ujt\|_{L^p} \leq Ce^{-c \no t 2^{2l}} \|\D_l \D_j u_0\|_{L^p}\\
+C_F \max(\nu, \nu') 2^{3j-l} e^{CV(t)}(e^{CV(t)}-1) \int_0^t e^{-c \no(t-\tau) 2^{2l}} \|u_j(\tau)\|_{L^p} d\tau\\
+C 2^{-|l-j|} e^{CV(\tau)} \Sum_{k\geq -1} 2^{-|k-j|} \int_o^t\|\n v(\tau)\|_{L^\infty} \|\D_k u(\tau)\|_{L^p} d\tau\\
+C \int_0^t e^{-c \no (t-\tau) 2^{2l}} \big(e^{CV(t)} 2^{-|l-j|} \|\Fj(\tau)\|_{L^p}+\|\Gj(\tau)\|_{L^p}\big) d\tau \bigg).
\end{multline}
Taking the $L_t^r$-norm and then summing for $l\geq j-N_0$, we obtain what follows (we skip the details as the method is exactly the same as in \cite{TH1}): for all $j\geq N_0$ and $t\in[0,T]$:
\begin{multline}
(\no r 2^{2j})^\frac{1}{r} \Sum_{l\geq j-N_0} \|\D_l \ujt\|_{L_t^r L^p} \leq C \|\D_j u_0\|_{L^p} +C_r e^{CV(t)}2^{N_0(1+\frac{2}{r})} \|\Fe\|_{L_t^1 L^p}\\
+\frac{C}{\no} 2^{2N_0} 2^{-j(2-\frac{2}{r})}\|\Ge\|_{L_t^r L^p} +C_F\frac{\max(\nu, \nu')}{\no} 2^{3N_0} e^{CV(t)}(e^{CV(t)}-1) (\no r 2^{2j})^\frac{1}{r} \|u_j\|_{L_t^r L^p}\\
+Ce^{CV(t)} \|\n v\|_{L_t^{\bar{r}} L^\infty} t^\frac{1}{r} 2^{N_0} \Sum_{k\geq -1} \no^\frac{1}{r} 2^\frac{2j}{r} 2^{-|k-j|} \|\D_k u\|_{L_t^r L^p}.
\label{estimHFbis}
\end{multline}
As in \eqref{estimBFR}, we need to estimate the low frequencies:
$$
\|S_{j-N_0} \ujt\|_{L_t^r L^p} \leq C e^{CV(t)} 2^{-N_0} \|\uj\|_{L_t^r L^p}.
$$
Multiplying by $2^{js}(\no r)^\frac{1}{r}$ and introducing $\aa(r,s)=\min\left(1+(s+\frac{2}{r}), 1-(s+\frac{2}{r})\right)>0$ (this is the reason of the restrictions  on the indexes) and:
$$
U_j(t) =\no^\frac{1}{r} 2^{j(s+\frac{2}{r})} \|\uj\|_{L_t^r L^p},
$$
we obtain that for all $j\geq N_0$,
\begin{multline}
U_j(t) \leq C 2^{js} \|\D_j u_0\|_{L^p} + \left( C2^{-N_0} e^{CV(t)} +C_F\frac{\max(\nu, \nu')}{\no} 2^{3N_0} e^{CV(t)}(e^{CV(t)}-1) \right) U_j(t)\\
+C e^{CV(t)} t^\frac{1}{r} \|\n v\|_{L_t^{\bar{r}} L^\infty} 2^{N_0} \sum_{k \geq -1} U_k(t) 2^{-\aa (r,s) |k-j|}\\
+C e^{CV(t)} 2^{3N_0} 2^{js} \|\Fj\|_{L_t^1 L^p} +\frac{C}{\no} 2^{2N_0} 2^{s+\frac{2}{r}-2} \|\Gj\|_{L_t^r L^p}.
\end{multline}
We can bound the low frequencies ($j\leq N_0$) as in \cite{TH1} and we finally end up with:
\begin{multline}
U_j(t) \leq C 2^{js} \|\D_j u_0\|_{L^p}(1+(\no t)^\frac{1}{r} 2^{N_0})\\
+\left( C2^{-N_0} e^{CV(t)} +C_F\frac{\max(\nu, \nu')}{\no} 2^{3N_0} e^{CV(t)}(e^{CV(t)}-1) \right) U_j(t)\\
+C \left(e^{CV(t)} 2^{3N_0} +(\no t)^\frac{1}{r} 2^{2N_0}\right) 2^{js}\|\Fj\|_{L_t^1 L^p}\\
+\left(\frac{C}{\no} 2^{2N_0} +\no^\frac{1}{r} t 2^{N_0}\right) 2^{s+\frac{2}{r}-2} \|\Gj\|_{L_t^r L^p}\\
+C \left( e^{CV(t)} 2^{N_0} t^\frac{1}{r} +(\no t)^\frac{1}{r} \right)\|\n v\|_{L_t^{\bar{r}} L^\infty} \sum_{k \geq -1} U_k(t) 2^{-\aa (r,s) |k-j|}.
\end{multline}
As there exists a constant $C_{r,s}>0$ so that we have:
$$
\sum_{k \geq -1} U_k(t) 2^{-\aa (r,s) |k-j|} \leq C_{r,s} \sup_{k \geq -1} U_k(t),
$$
taking the supremum for $j\geq -1$, we choose $T$ small enough and $N_0$ large enough so that:
$$
\begin{cases}
\vspace{1mm}
\displaystyle{CC'T^\frac{1}{4} \leq \frac{1}{2}\no^\frac{3}{4},}\\
\vspace{1mm}
\displaystyle{e^{CV(T)}-1 \leq \frac{1}{2},}\\
\vspace{1mm}
\displaystyle{C_F\frac{\max(\nu, \nu')}{\min(\nu, \nu')} 2^{3N_0} e^{CV(T)} (e^{CV(T)}-1) \leq \frac{1}{6},}\\
\vspace{1mm}
\displaystyle{C2^{-N_0} e^{CV(T)}\leq \frac{1}{6},}\\
\displaystyle{C_{r,s} \left( e^{CV(T)} 2^{N_0} +\no^\frac{1}{r} \right) T^\frac{1}{r} \|\n v\|_{L_T^{\bar{r}} L^\infty} \leq \frac{1}{6}.}
\end{cases}
$$ 
These conditions are realized when we fix $N_0 \in\N$ so large that $2^{N_0} \geq 9C$ and $T$ is chosen small enough so that:
$$
\begin{cases}
\vspace{1mm}
\displaystyle{T \leq \frac{\no^3}{(2CC')^4},}\\
\vspace{1mm}
\displaystyle{e^{CV(T)}-1 \leq \min \left(\frac{1}{2}, \frac{1}{9 C_F \frac{\max(\nu, \nu')}{\min(\nu, \nu')} 2^{3N_0}}\right),}\\
\displaystyle{T+\int_0^T \|\n v\|_{L^\infty}^{\bar{r}} d\tau \leq \frac{C_{r,s}}{6\left(2^{N_0+1} +\no^\frac{1}{r} \right)}.}
\end{cases}
$$ 
Under these conditions we obtain that for all $t\in[0,T]$,
$$
\sup_{j \geq -1} U_j(t) \leq C_{r,\no, F} \left(\|u_0\|_{B_{p,\infty}^s} +(1+(\no t)^\frac{1}{r}) \|\Fe\|_{\tilde{L}_t ^1 B_{p,\infty}^s} +C(\frac{1}{\no} +\no^\frac{1}{r} t)\|\Ge\|_{\tilde{L}_t ^r B_{p,\infty}^{s+\frac{2}{r}-2}}\right),
$$
which leads to the conclusion of Theorem \ref{thCs}.

\section{Appendix}

The first part is devoted to a quick presentation of the Littlewood-Paley theory. In the second section we briefly recall general considerations on flows. The last section gives results on the diffeomorphisms introduced in the proofs of the theorems.

\subsection{Besov spaces}

As usual, the Fourier transform of $u$ with respect to the space variable is denoted by $\mathcal{F}(u)$ or $\hat{u}$. 
In this section we will state classical definitions and properties concerning the homogeneous and nonhomogeneous dyadic decomposition with respect to the Fourier variable. For a complete presentation we refer to \cite{Dbook} (Chapter 2).

Let $\chi$ a smooth radial function supported in the ball $B(0,\frac{4}{3})$, equal to 1 in a neighborhood of $B(0,\frac{3}{4})$ and such that $r\mapsto \chi(r.e_r)$ is nonincreasing over $\R_+$. If we set $\varphi(\xi)=\chi(\xi/2)-\chi(\xi)$, then $\varphi$ is compactly supported in the annulus $\cC=\{\xi\in \R^d, c_0=\frac{3}{4}\leq |\xi|\leq C_0=\frac{8}{3}\}$ and we also have that:
\begin{equation}
\begin{cases} 
 \forall \xi\in \R^d\setminus\{0\}, \quad \displaystyle{\sum_{l\in\Z} \varphi(2^{-l}\xi)=1},\\
 \forall \xi\in \R^d, \quad \chi(\xi) +\displaystyle{\sum_{l \geq 0} \varphi(2^{-l}\xi)=1}.
\end{cases}
\label{LPxi}
\end{equation}
From this we define the homogeneous dyadic blocks: for all $j\in \Z$,
$$
\ddj u:= \varphi(2^{-j}D)u =2^{jd} h(2^j.)* u, \quad \mbox{with } h=\cF^{-1} \varphi.
$$
We recall that $\hat{\phi(D)u}(\xi)=\phi(\xi) \hat{u} (\xi)$. Then we have for tempered distributions (modulo polynomials, we refer to \cite{Dbook} for precisions):
\begin{equation}
u=\Sum_l \ddl u
\label{LPsomme} 
\end{equation}
Similarly, the nonhomogeneous dyadic operators are defined for all $j \in \Z$, by
$$
\begin{cases}
\displaystyle{\forall j\leq -2, \quad \D_j u=0,}\\
\displaystyle{\D_{-1}u =\chi(D)u,}\\
\displaystyle{\forall j\geq 0, \quad \D_j u:= \varphi(2^{-j}D)u =2^{jd} h(2^j.)* u.}
\end{cases}
$$
We also define the low frequency cut-off operators
$$
\begin{cases}
\dot{S}_j u:= \chi(2^{-j}D) u=\Sum_{q\leq j-1} \ddq u =2^{jd} k(2^j.)* u$ if $k=\cF^{-1} \chi,\\
S_j u=\chi(2^{-j}D)u =\Sum_{q\leq j-1} \D_q u.
\end{cases}
$$
We can now define the nonhomogeneous Besov spaces:
\begin{defi}
\sl{For $s\in\R$ and $1\leq p,r\leq\infty,$ we set
$$
\|u\|_{B^s_{p,r}}:=\bigg(\sum_{l\geq -1} 2^{rls}
\|\D_l u\|^r_{L^p}\bigg)^{\frac{1}{r}}\ \text{ if }\ r<\infty
\quad\text{and}\quad
\|u\|_{B^s_{p,\infty}}:=\sup_{l} 2^{ls}
\|\D_l u\|_{L^p}.
$$
The nonhomogeneous Besov space $B^s_{p,r}$ is the subset of tempered distributions such that $\|u\|_{B^s_{p,r}}$ is finite.
}
\end{defi}
It is more delicate for the homogeneous spaces:
\begin{defi}
\sl{For $s\in\R$ and $1\leq p,r\leq\infty,$ we set
$$
\|u\|_{\dot B^s_{p,r}}:=\bigg(\sum_{l\in \Z} 2^{rls}
\|\ddl u\|^r_{L^p}\bigg)^{\frac{1}{r}}\ \text{ if }\ r<\infty
\quad\text{and}\quad
\|u\|_{\dot B^s_{p,\infty}}:=\sup_{l} 2^{ls}
\|\ddl u\|_{L^p}.
$$
The homogeneous Besov space $\dot B^s_{p,r}$ is the subset of tempered distributions such that $\lim_{j \rightarrow -\infty} \|\dot{S}_j u\|_{L^\infty}=0$ and $\|u\|_{\dot B^s_{p,r}}$ is finite.
}
\end{defi}
Once more, we refer to \cite{Dbook} (chapter $2$) for properties of the nonhomogeneous and homogeneous Besov spaces. For now let us just state that:
\begin{itemize}
\item there exists a constant $C$ such that for all $j\in \Z$, $p\in[1,\infty]$ and all $u$,
$$
\|\D_j u\|_{L^p} +\|\ddj u\|_{L^p} +\|S_j u\|_{L^p} +\|\dot{S}_j u\|_{L^p} \leq C\|u\|_{L^p}.
$$
\item for all $j,j'$ such that $|j-j'|\geq 2$, $\ddj \dot{\D}_{j'} =0$ (the same for nonhomogeneous operators).
\end{itemize}
The Littlewood-Paley decomposition allows us to work with spectrally localized (therefore smooth) functions rather than with rough objects. We obtain bounds for each dyadic block in spaces of type $L^\rho_T L^p$. Then, with a view to get estimates in  $L^\rho_T \dot B^s_{p,r}$, we perform a summation in $\ell^r(\Z)$. In fact we do not bound the $L^\rho_T \dot B^s_{p,r}$ norm as the time integration has been performed before the $\ell^r$ summation.
This leads to the definition of the spaces $\tilde L^\rho_T \dot B^s_{p,r}$ from the following norm:
\begin{defi}\label{d:espacestilde}
For $T>0,$ $s\in\R$ and  $1\leq r,\rho\leq\infty,$
 we set
$$
\|u\|_{\tilde L_T^\rho \dot B^s_{p,r}}:=
\bigl\Vert2^{js}\|\ddq u\|_{L_T^\rho L^p}\bigr\Vert_{\ell^r(\Z)}.
$$
\end{defi}
The spaces $\tilde L^\rho_T \dot B^s_{p,r}$ can be compared with the spaces  $L_T^\rho \dot B^s_{p,r}$ thanks to the Minkowski inequality: we have
$$
\|u\|_{\tilde L_T^\rho \dot B^s_{p,r}}
\leq\|u\|_{L_T^\rho \dot B^s_{p,r}}\ \text{ if }\ r\geq\rho\quad\hbox{and}\quad
\|u\|_{\tilde L_T^\rho \dot B^s_{p,r}}\geq
\|u\|_{L_T^\rho \dot B^s_{p,r}}\ \text{ if }\ r\leq\rho.
$$
All the properties of continuity for the product and composition which are true in Besov spaces remain true in the above  spaces.
\medbreak
Let us now recall a few nonlinear estimates in Besov spaces. Formally, any product of two distributions $u$ and $v$ may be decomposed into 
\begin{equation}\label{eq:bony}
uv=T_uv+T_vu+R(u,v), \mbox{ where}
\end{equation}
$$
T_uv:=\sum_l\dot S_{l-1}u\ddl v,\quad
T_vu:=\sum_l \dot S_{l-1}v\ddl u\ \hbox{ and }\ 
R(u,v):=\sum_l\sum_{|l'-l|\leq1}\ddl u\,\dot\Delta_{l'}v.
$$
The above operator $T$ is called a ``paraproduct'' whereas $R$ is called a ``remainder''. The decomposition \eqref{eq:bony} has been introduced by J.-M. Bony in \cite{Bony}. We refer to \cite{Dbook} for general properties and for paraproduct and remainder estimates. The same can be defined for nonhomogeneous spaces.

\subsection{Estimates for the flow of a smooth vector-field}
In this section, we recall classical estimates for the flow associated to the smooth vectorfield $S_{j-1} v$. We refer to \cite{Dbook} (sections $2.1.3$ and $3.1.2$), \cite{CD, Dlagrangien, CHVP, Corder} for more details in the compressible case. We also refer to \cite{TH1} for the incompressible Navier-Stokes case.
\begin{prop}
\label{p:flow}
\sl{
Under the assumptions of theorems \ref{thC2} and \ref{thCs}, if $V(t) :=\int_0^t\|\nabla v(t')\|_{L^\infty}\,dt'.$  Let $\psj$ the associated flow:
$$
\psj(x)=x+\int_0^t S_{j-1}v(t',\psjto(x))\,d\tau'.
$$
Then for all $t\in\R,$ the flow $\psj$ is a smooth volume-preserving diffeomorphism over $\R^3$ and there exists a constant $C>0$ such that one has, if $t\geq 0$,
$$
\begin{cases}
\vspace{1mm}
\|D\psi_t^{\pm1}\|_{L^\infty}\leq e^{CV(t)},\\
\vspace{1mm}
\|D\psi_t^{\pm1}-I_d\|_{L^\infty}\leq  e^{CV(t)}-1,\\
\forall k\geq 2\quad \|D^k\psi_t^{\pm1}\|_{L^\infty}\leq C2^{(k-1)j}\big(e^{CV(t)}-1\big),
\end{cases}
$$
}
\end{prop}
where $Df$ denotes the jacobian matrix of $f$.

\subsection{Study of the diffeomorphisms $m_x(\pm.)$}
In this section we prove precise results and estimates for the diffeomorphisms $m_x(\pm.)$ and for the quantities $Y_\pm$ defined as follows (see \eqref{notmxY}):
$$
\begin{cases}
\vspace{1mm}
m_x(y)=\psjm(x)-\psjm(x+y),\\
\displaystyle{Y_{\pm}=\frac{|m_x(\pm y)|}{|y|} \quad \mbox{and}\quad Y_{\pm}^F=\frac{|m_x(\pm y)|_\frac{1}{F}}{|y|_\frac{1}{F}} \quad \mbox{if } y\neq 0}.
\end{cases}
$$
First, let us state the following result:
\begin{prop}
\sl{There exists a constant $C>0$ such that for all $x,y\in\R^3$,
\begin{enumerate}
\item $e^{-C V(t)} \leq Y_{\pm} \leq e^{C V(t)}$, where $V(t)=\int_0^t \|\n v(\tau)\|_{L^\infty} d\tau$,
\item or equivalently, $e^{-C V(t)}|y| \leq |m_x(\pm y)| \leq e^{C V(t)}|y|$,
\item $|Y_\pm-1| \leq e^{2C V(t)}-1$ and $|\frac{1}{Y_\pm}-1| \leq e^{2C V(t)}-1$,
\item or equivalently, $\bigg||m_x(\pm y)|-|y|\bigg| \leq (e^{2C V(t)}-1) |y|$.
\item $|Y_+-Y_-| \leq e^{2C V(t)}(e^{2C V(t)}-1) \min(1,2^j|y|),$
\item or equivalently, $\bigg||m_x(-y)|-|m_x(y)|\bigg| \leq  e^{2C V(t)}(e^{2C V(t)}-1) \min(1,2^j|y|)|y|.$
\item The same for $Y_{\pm}^F$ (where the $\R^3$ euclidiean norm $|.|$ is replaced by $|.|_\frac{1}{F}$, see \eqref{normodif} and we).
\item $\bigg|m_x(\pm y)\pm y|\bigg| \leq e^{2C V(t)}(e^{2C V(t)}-1) |y|$.
\item $\bigg|m_x(-y)+m_x(y)|\bigg| \leq  e^{2C V(t)}(e^{2C V(t)}-1) \min(1,2^j|y|)|y|.$
\end{enumerate}
}
\label{propmx}
\end{prop}
\begin{rem}
\sl{We emphasize that as $\psj$ is close to $Id$, $m_x(-y)\sim y$ and $m_x(y)\sim -y$.}
\end{rem}
\textbf{Proof :} We refer to \cite{CHVP, Corder} for the proof of points $1$ to $7$. We will prove here the last two points. Let us begin with point $8$ in the case of $+y$ (the case $-y$ is dealt exactly the same way): thanks to the definition of $\psj$ (see \eqref{defpsi}), for all $y\in \R^3$,
$$
\psj (x)-x =\int_0^t S_{j-1}v\big(\tau, \psjto (x)\big)d\tau.
$$
Applying this at the point $\psjm(x)$, we obtain:
$$
x-\psjm (x) =\int_0^t S_{j-1}v\big(\tau, \psjto \circ \psjm (x)\big)d\tau.
$$
Now, applying this to the point $x-y$ and thanks to the definition of $m_x$ (see \eqref{notmxY}0, we end up with:
\begin{multline}
y-m_x(-y)= y- \big(\psjm (x)- \psjm(x-y)\big) =-\left( \psjm (x)-x-\big(x-y-\psjm(x-y)\big)\right)\\
=-\int_0^t \left(S_{j-1}v\big(\tau, \psjto \circ \psjm (x)\big) -S_{j-1}v\big(\tau, \psjto \circ \psjm (x-y)\big) \right)d\tau.
\end{multline}
Thanks to the Taylor formula, the integrand is equal to:
$$
\int_0^1 \n S_{j-1}v\big(\tau, \psjto \circ \psjm (x-y+sy)\big)\cdot \n \psjto \circ \psjm (x-y+sy) \cdot \n \psjm(x-y+sy) \cdot y ds,
$$
which allows us to obtain (using Proposition \ref{p:flow}):
\begin{multline}
|y-m_x(-y)| \leq \int_0^t \|\n S_{j-1}v\|_{L^\infty} \|\n \psjto\|_{L^\infty} \|\n \psjm\|_{L^\infty} |y| d\tau\\
\leq V(t) e^{2C V(t)} |y| \leq e^{2C V(t)}(e^{C V(t)}-1) |y|.
\end{multline}
which end the proof of point $8$. For the last point we simply write that:
$$
m_x(-y) +m_x(y)= 2\psjm (x) -\psjm(x-y) -\psjm(x+y).\\
$$
As in \cite{CHVP, Corder}, an elementary use of the Taylor formula implies that for a smooth $\cC^2$ function,
$$
\sup_{x\in \R^3} |2f(x)-f(x-y)-f(x+y)| \leq |y|^2 \|\n^2 f\|_{L^\infty},
$$
and then, thanks to Proposition \ref{p:flow}, we obtain that for all $x,y\in R^3$ with $y\neq 0$,
$$
|m_x(-y) +m_x(y)| \leq |y|^2 \|\n^2 \psjm\|_{L^\infty} \leq e^{CV(t)} \big(e^{CV(t)} -1\big) 2^j |y|^2.
$$
On the other hand, thanks to point $8$ from the proposition,
\begin{multline}
|m_x(-y) +m_x(y)|= |m_x(-y)-y+ y+m_x(y)| \leq |m_x(-y)-y| +|m_x(y)+y|\\
\leq 2e^{2CV(t)} \big(e^{2CV(t)} -1\big) |y|.
\end{multline}
Gathering the last two estimates ends the proof of point $9$. $\blacksquare$
\\

We will conclude this section with the following result:
\begin{prop}
\sl{With the same notations as before, if $t$ is so small that:
$$
e^{2C V(t)}-1 \leq \frac{1}{2},
$$
then for all smooth function $k$ defined on $[\frac{1}{2}, \frac{3}{2}]$, and all $x,y \in \R^3$ with $y\neq 0$,
$$
\begin{cases}
\vspace{1mm}
\displaystyle{|k(1) -k(Y_\pm^{(F)})| \leq \|k'\|_{L^\infty([\frac{1}{2}, \frac{3}{2}])} \big(e^{2C V(t)}-1\big),}\\
\displaystyle{|k(Y_-^{(F)}) -k(Y_+^{(F)})| \leq \|k'\|_{L^\infty([\frac{1}{2}, \frac{3}{2}])} e^{2C V(t)}(e^{2C V(t)}-1) \min(1,2^j|y|).}
\end{cases}
$$
}
\label{propkY}
\end{prop}
\textbf{Proof :} an elemental use of the Taylor formula gives that:
$$
k(1) -k(Y_\pm^{(F)}) =\int_0^1 k'\big(1+(1-s)(Y_\pm^{(F)}-1)\big)\cdot (1-Y_\pm^{(F)}) ds,
$$
so immediately, we obtain
$$
|k(1) -k(Y_\pm^{(F)})| \leq \|k'\|_{L^\infty([1-(e^{2C V(t)}-1), 1+(e^{2C V(t)}-1)])} |1-Y_\pm^{(F)}|.
$$
Similarly, we obtain that:
$$
|k(Y_-^{(F)}) -k(Y_+^{(F)})| \leq \|k'\|_{L^\infty([1-(e^{2C V(t)}-1), 1+(e^{2C V(t)}-1)])} |Y_-^{(F)}-Y_+^{(F)}|
$$
and the conclusion follows from the assumption on $t$ and Proposition \ref{propmx}. $\blacksquare$

\subsection{Results on the kernel $K$}

In this section we wish to gather various results on the kernel involved in the singular integral form of the QG operator $\G$ (see \eqref{Kexpr}). We begin by rewriting $K$ as follows: $\forall y\in \R^3$
$$
K(y)= -\frac{2C}{F^3}\frac{y_1^2+ y_2^2 -\frac{3}{F^2}y_3^2}{\big(y_1^2+ y_2^2 +\frac{1}{F^2}y_3^2\big)^3} =C_F \frac{\big(y|L(y)\big)}{|y|_{\frac{1}{F}}^6},
$$
where $(.|.)$ is the cannonical inner product on $\R^3$ and $L$ is the isomorphism of $\R^3$ defined for all $y\in \R^3$ by
$$
L(y_1, y_2, y_3)=(y_1, y_2, -\frac{3}{F^2}y_3).
$$
First we have the following estimates:
$$
\begin{cases}
\min(1, \frac{3}{F^2})|y| \leq |L(y)| \leq \max(1, \frac{3}{F^2}) |y|\\
\min(1, \frac{3}{F})|y|_\frac{1}{F} \leq |L(y)| \leq \max(1, \frac{3}{F}) |y|_\frac{1}{F}\\
\min(1, \frac{1}{F})|y| \leq |y|_\frac{1}{F} \leq \max(1, \frac{1}{F}) |y|.
\end{cases}
$$
As $L$ is linear we directly write:
\begin{lem}
\sl{For all $a,b\in \R^3$, we have:
\begin{equation}
\big(b|L(b)\big)-\big(a|L(a)\big) =\big(b+a|L(b-a)\big) =\big(b-a|L(b+a)\big).
\label{Lab}
\end{equation}
}
\end{lem}
Contrary to the cases of \cite{CHVP, Corder} the kernel has no sign and reaches the value $0$, so we will have to be much more careful than in these works.
\begin{prop}
\sl{With the same notations as in Proposition \ref{propmx}, there exists a constant $C_F$ such that for all $t$ so small that
$$
e^{2C V(t)}-1 \leq \frac{1}{2},
$$
then we have:
\begin{equation}
\begin{cases}
\vspace{1mm}
\displaystyle{|K(y)-K\big(m_x(-y)\big)| \leq \frac{C_F}{|y|^4} e^{6CV(t)} \big(e^{2CV(t)}-1\big),}\\
\displaystyle{|K\big(m_x(y)\big)-K\big(m_x(-y)\big)| \leq \frac{C_F}{|y|^4} e^{6CV(t)} \big(e^{2CV(t)}-1\big) \min(1, 2^j|y|).}
\end{cases}
\end{equation}
}
\label{propKymx}
\end{prop}

\textbf{Proof :} Using \eqref{Lab}, we immediately obtain that for any $y,z \in \R^3\setminus \{0\}$,
\begin{multline}
K(y)-K(z)=C_F \left(\frac{\big(y|L(y)\big)}{|y|_{\frac{1}{F}}^6} -\frac{\big(z|L(z)\big)}{|z|_{\frac{1}{F}}^6}\right)\\
=C_F \left(\frac{1}{|y|_{\frac{1}{F}}^6}\big(y+z|L(y-z)\big) +\big(z|L(z)\big) \big(\frac{1}{|y|_{\frac{1}{F}}^6} -\frac{1}{|z|_{\frac{1}{F}}^6}\big)\right).
\end{multline}
This implies that for all $y\in \R^3\setminus \{0\}$,
\begin{multline}
K(y)-K\big(m_x(-y)\big)\\
=\frac{C_F}{|y|_{\frac{1}{F}}^6} \left(\big(y+m_x(-y)|L(y-m_x(-y))\big) +\big(m_x(-y)|L(m_x(-y))\big) \big(1-\frac{1}{(Y_-^F)^6}\big)\right).
\end{multline}
so that
\begin{multline}
|K(y)-K\big(m_x(-y)\big)|\\
\leq\frac{C_F}{|y|^6} \left(|y+m_x(-y)|\cdot C_F|y-m_x(-y)| + C_F |m_x(-y)|^2 \big(1-\frac{1}{(Y_-^F)^6}\big)\right).
\end{multline}
and similarly
\begin{multline}
|K\big(m_x(-y)\big)-K\big(m_x(y)\big)| \leq\frac{C_F}{|m_x(-y)|_{\frac{1}{F}}^6} |m_x(-y)+m_x(y)|\cdot |m_x(-y)-m_x(y)|\\
+C_F \frac{|m_x(y)|^2}{|y|_{\frac{1}{F}}^6} \left(\frac{1}{(Y_-^F)^6} -\frac{1}{(Y_+^F)^6}\big)\right).
\end{multline}
Finally, the conclusion is obtained thanks to Propositions \ref{propmx} and \ref{propkY} (with $k(u)=u^{-6}$). $\blacksquare$
\\
\\
The author wishes to thank T. Hmidi, C. Imbert and G. Karch for useful discussions.

\end{document}